
\magnification=1200

\font\tengoth=eufm10
\font\sevengoth=eufm7
\newfam\gothfam
\textfont\gothfam=\tengoth
\scriptfont\gothfam=\sevengoth
\def\goth{\fam\gothfam\tengoth}

\font\tenbboard=msbm10
\font\sevenbboard=msbm7
\newfam\bboardfam
\textfont\bboardfam=\tenbboard
\scriptfont\bboardfam=\sevenbboard
\def\bboard{\fam\bboardfam\tenbboard}

\newif\ifpagetitre
\newtoks\auteurcourant \auteurcourant={\hfill}
\newtoks\titrecourant \titrecourant={\hfill}

\pretolerance=500 \tolerance=1000 \brokenpenalty=5000
\newdimen\hmargehaute \hmargehaute=0cm
\newdimen\lpage \lpage=14.3cm
\newdimen\hpage \hpage=20cm
\newdimen\lmargeext \lmargeext=1cm
\hsize=11.25cm
\vsize=18cm
\parskip=0cm
\parindent=12pt

\def\margehaute{\vbox to \hmargehaute{\vss}}
\def\margebasse{\vss}

\output{\shipout\vbox to \hpage{\margehaute\nointerlineskip
  \corpsdepage\margebasse}
  \advancepageno \global\pagetitrefalse
  \ifnum\outputpenalty>-2000 \else\dosupereject\fi}

\def\corpsdepage{\hbox to \lpage{\hss\pagetexte\hskip\lmargeext}}
\def\pagetexte{\vbox{\makeheadline\pagebody\makefootline}}
\headline={\ifpagetitre\titleheadline \else
  \ifodd\pageno\rightheadline \else \leftheadline\fi\fi}
\def\leftheadline{\hfil\the\auteurcourant\hfil}
\def\rightheadline{\hfil\the\titrecourant\hfil}
\def\titleheadline{\hfill}
\pagetitretrue

\font\petcap=cmcsc10
\def\pc#1#2|{{\tenrm#1\sevenrm#2}}
\def\pd#1 {{\pc#1|}\ }

\def\pointir{\discretionary{.}{}{.\kern.35em---\kern.7em}\nobreak\hskip 0em
 plus .3em minus .4em }

\def\titre#1|{\message{#1}
             \par\vskip 30pt plus 24pt minus 3pt\penalty -1000
             \vskip 0pt plus -24pt minus 3pt\penalty -1000
             \centerline{\bf #1}
             \vskip 5pt
             \penalty 10000 }
\def\section#1|
               {\par\vskip .3cm
               {\bf #1}\pointir}
\def\ssection#1|
             {\par\vskip .2cm
             {\it #1}\pointir}

\long\def\th#1|#2\finth{\par\vskip 5pt
              {\petcap #1\pointir}{\it #2}\par\vskip 3pt}
\long\def\tha#1|#2\fintha{\par\vskip 5pt
               {\petcap #1.}\par\nobreak{\it #2}\par\vskip 3pt}

\def\rem#1|
{\par\vskip 5pt
                {{\it #1}\pointir}}
\def\rema#1|
{\par\vskip 5pt
                {{\it #1.}\par\nobreak }}

\def\proof{\noindent {\it Proof. }}

\def\qed{\quad\raise -2pt \hbox{\vrule\vbox to 10pt
{\hrule width 4pt
\vfill\hrule}\vrule}}

\def\cqfd{\ifmmode
\unkern\quad\hfill\qed
\else
\unskip\quad\hfill\qed\bigskip
\fi}

\newcount\n
\def\exo{\advance\n by 1 \par \vskip .3cm {\bf \the \n}. }

\def\Det{\mathop{\rm Det}\nolimits}

\def\tr{\mathop{\rm tr}\nolimits}

\def\Tr{\mathop{\rm Tr}\nolimits}

\def\ad{\mathop{\rm ad}\nolimits}

\def\Ad{\mathop{\rm Ad}\nolimits}

\def\demi{{\textstyle{1\over 2}}\,}

\font\msbmten=msbm10
\def\ltimes{\mathbin{\hbox{\msbmten\char'156}}}


\hfill 14 January 2011

\vskip 20 pt

\centerline{ANALYSIS OF THE BRYLINSKI-KOSTANT MODEL}

\vskip 2pt

\centerline{FOR SPHERICAL MINIMAL REPRESENTATIONS}

\vskip 20 pt
 
\centerline{Dehbia ACHAB  \&  Jacques FARAUT}

\vskip 15pt

\noindent
{\bf Abstract} 
{\it We revisit with another view point the construction by R. Brylinski and B. Kostant of minimal representations of simple Lie groups. We start from a pair $(V,Q)$, where $V$ is a complex vector space and $Q$ a homogeneous polynomial of degree 4 on $V$.
The manifold $\Xi $ is an orbit of a covering of ${\rm Conf}(V,Q)$, the conformal group of the pair $(V,Q)$, in a finite dimensional representation space.  
By a generalized Kantor-Koecher-Tits construction we obtain a complex simple Lie algebra $\goth g$, and furthermore a real
form ${\goth g}_{\bboard R}$. The connected and simply connected Lie group $G_{\bboard R}$ with ${\rm Lie}(G_{\bboard R})={\goth g}_{\bboard R}$ acts unitarily on a Hilbert space of holomorphic functions defined on the manifold $\Xi $.}

\bigskip


\bigskip

\noindent
{\it Mathematics Subject Index 2010}:17C36, 22E46, 32M15, 33C80.

\bigskip

Introduction

1. The conformal group and the representation $\kappa $

2. The orbit $\Xi $ and the irreducible $K$-invariant Hilbert \hfill \break subspaces of ${\cal O}(\Xi )$

3. Decomposition into simple Jordan algebras

4. Generalized Kantor-Koecher-Tits construction

5. Representation of the generalized Kantor-Koecher-Tits Lie algebra

6. The unitary representation of the Kantor-Koecher-Tits group

References

\vskip 1cm

\section Introduction|
The construction of a realization for the minimal unitary representation of a simple Lie group by using geometric quantization has been the topic of many papers during the last thirty years: [Rawnsley-Sternberg,1982], [Torasso,1983], and more recently [Kobayashi-\O rsted,2003].
In  a series of papers R. Brylinski and B. Kostant have
introduced and studied  a geometric quantization  of minimal nilpotent orbits for
simple real Lie groups which are not of Hermitian type:
[Brylinski-Kostant,1994,1995,1997], [Brylinski, 1997,1998] . They have constructed the
associated  irreducible unitary  representation on a Hilbert space of 
half forms on
the minimal nilpotent orbit. This can be considered as a Fock model for the minimal representation.
In this paper we revisit this construction with another point of view. We
start from a pair $(V,Q)$ where $V$ is a complex vector space  and $Q$ is a
homogeneous polynomial on $V$ of degree $4$. The structure group  Str$(V,Q)$, 
for which
$Q$ is a semi-invariant, is assumed to have a symmetric open orbit. The
conformal group Conf$(V,Q)$ consists of rational transformations of $V$ whose
differential belongs to Str$(V,Q)$.
The main geometric object is the orbit $\Xi$ of $Q$ under $K$, a covering
of Conf$(V,Q)$, on a space $\cal W$ of polynomials on $V$.
Then, by a generalized Kantor-Koecher-Tits construction, starting from
the Lie algebra $\goth k$ of $K$, we obtain a simple Lie algebra $\goth g$
such that the pair $(\goth g,\goth k)$ is non Hermitian.
As a vector space $\goth g=\goth k \oplus \goth p$, with $\goth p=\cal W$. The
main point is to define a bracket 
$${\goth p}\oplus{\goth p}\rightarrow {\goth k},\quad (X,Y)\mapsto [X,Y],$$ 
such that $\goth g$ becomes a Lie algebra. The Lie algebra
$\goth g$ is 5-graded:
$$\goth g=\goth g_{-2}\oplus\goth g_{-1}\oplus \goth g_0 \oplus \goth g_1\oplus\goth g_2.$$
In the fourth part one defines a representation  $\rho$  of $\goth g$ on the
space ${\cal O}(\Xi)_{\rm fin}$ of polynomial functions on $\Xi$. In a first step one
defines a representation of an ${\goth sl}_2$-triple $(E,F,H)$. It turns out that this
is only possible under a condition $T$. In such a case one obtains an
irreducible unitary representation of the connected and simply connected
group $\widetilde G_{\bboard R}$ whose Lie algebra is a real form of $\goth g$. The
representation is spherical.
 It is realized on a Hilbert space of holomorphic functions on
$\Xi  $. There is an explicit formula for the  reproducing kernel of $\cal H$
involving a hypergeometric function ${}_1F_2$.
Further the space $\cal H$ is a weighted Bergman space with a weight
taking in general both positive and negative values.
\vskip 1pt
If  $Q=R^2$ or $Q=R^4$ where $R$ is a semi-invariant, then by considering a
covering of order 2 or 4  of the orbit $\Xi$, one can obtain one or 3 other
unitary representations of $\widetilde G_{\bboard R}$. They are not spherical.
If the condition T is not satisfied, by a modified construction, one
still obtains an irreducible representation of $\widetilde G_{\bboard R}$ which is
not spherical. This last point is the subject of a paper in preparation by the first author.
\vskip 1pt
The construction of a Schr\"odinger model for the minimal representation of the group $O(p,q)$ is the subject of a recent book by T. Kobayashi and G. Mano [2008]. We should not wonder that there is a link between both models: the Fock and the Schr\"odinger models, and that there is an analogue of the Bargmann transform in this setting.

\section 1. The conformal group and the representation $\kappa $|
Let $V$ be a finite dimensional complex vector space and $Q$ a 
homogeneous polynomial on $V$. Define
$$L={\rm Str}(V,Q)=\{g\in GL(V) \mid \exists \gamma=\gamma(g),
Q(g\cdot x)=\gamma(g)Q(x)\}.$$
Assume that there exists $e \in V$ such that 
\vskip 1pt
\qquad (1)  The symmetric bilinear form
$$\langle x,y \rangle =-D_xD_y\log Q(e),$$
is non-degenerate.
\vskip 1pt
\qquad (2) The orbit $\Omega=L\cdot e$ is open.
\vskip 1pt
\qquad (3) The orbit $\Omega=L\cdot e$ is symmetric, i.e.  the pair
$(L,L_0)$, with $L_0=\{g \in L \mid g\cdot e=e\},$
is symmetric, which means that there is an involutive 
automorphism $\nu $ of $L$ 
such that $L_0$ is open in $\{g\in L \mid \nu (g)=g\}$.

We will equip the vector space $V$ with a Jordan algebra structure.
The Lie algebra ${\goth  l}={\rm Lie}(L)$ of $L={\rm Str}(V,Q)$ decomposes into the $+1$ and $-1$ eigenspaces of the differential of $\nu $ :
${\goth  l}={\goth l}_0+{\goth q}$,
where ${\goth l}_0=\{X\in {\goth l} \mid X.e=e\}= {\rm Lie}(L_0)$. Since
the orbit $\Omega$ is open, the map
$${\goth q} \rightarrow V, \quad X\mapsto X.e,$$
is a linear isomorphism. If $X\cdot e=x \quad (X\in {\goth q}, x\in V)$ one
writes $X=T_x$. The product on $V$ is defined by
$$xy=T_x\cdot y=T_x\circ T_y\cdot e.$$

\th Theorem 1.1|
This product makes $V$ into a
semi-simple complex Jordan algebra:
\vskip 1pt
{\rm (J1)} For $x,y\in V, xy=yx$.
\vskip 1pt
{\rm (J2)} For $x,y\in V, x^2(xy)=x(x^2y)$.
\vskip 1pt
{\rm (J3)} The symmetric bilinear form $\langle.,.\rangle$ is
associative:
$$\langle xy,z\rangle =\langle x,yz\rangle .$$
\finth

\proof 
(a) This product is commutative. In fact 
$$xy-yx=[T_x,T_y]\cdot e=0,$$
since $[{\goth q},{\goth q}] \subset {\goth l}_0$.
\vskip 1pt
(b) Let $\tau$ be the differential of $\gamma$ at the identity
element of $L$: for $X \in {\goth l}$,
$$\tau(X)={d\over dt}\Big|_{t=0} \gamma(\exp tX).$$

\tha Lemma 1.2| 
$$\eqalign{
{\rm (i)} \quad & (D_x\log Q)(e)=\tau(T_x), \cr
{\rm (ii)} \quad & (D_xD_y\log Q)(e)=-\tau(T_{xy}),\cr
{\rm (iii)} \quad & (D_xD_yD_z\log Q)(e)=\demi \tau(T_{(xy)z}).\cr}$$
\fintha

The proof amounts to differentiating at $e$ the relation 
$$\log Q(\exp T_x\cdot e)=\tau(T_x)+\log Q(e),$$
up to third order. (See Exercise 5 in [Satake, 1980], p.38.)
Hence, by (ii), $\langle x,y\rangle=\tau(T_{xy})$, and, by (iii), the
symmetric bilinear form $\langle .,.\rangle$ is associative.
\vskip 2pt
\qquad (c) Define the associator of three elements $x,y,z$ in
$V$ by 
$$[x,y,z]=x(zy)-(xz)y=[L(x),L(y)]z.$$
Identity (J2) can be written: $[x^2,y,x]=0$ for all $x,y \in V$.
It can be shown by following the proof of Theorem 8.5 in
[Satake,1980], p.34, which is also the proof of Theorem III.3.1
in [Faraut-Koranyi,1994], p.50. 
\hfill \qed

\vskip 6pt

The Jordan algebra $V$ is a direct sum of simple ideals: 
$$V=\bigoplus_{i=1}^sV_i,$$
and
$$Q(x)=\prod _{i=1}^s\Delta_i(x_i)^{k_i}\quad
(x=(x_1,\ldots,x_s)),$$
where $\Delta_i$ is the
determinant polynomial of the simple Jordan algebra 
$V_i$ and the $k_i$ are positive integers. The degree of $Q$ is equal
to $\sum _{i=1}^sk_ir_i$, where $r_i$ is the rank of $V_i$.

\vskip 6pt

The conformal group ${\rm Conf}(V,Q)$ is the group of rational transformations $g$ of $V$ generated by: the translations $z\mapsto z+a$ ($a\in V$), the dilations $z\mapsto \ell \cdot z$
($\ell \in L$), and the inversion $j:z\mapsto -z^{-1}$. A transformation $g\in {\rm Conf}(V,Q)$
is conformal in the sense that the differential $Dg(z)$ belongs to $L\in {\rm Str}(V,Q)$ at any point $z$ where $g$ is defined.

Let $\cal W$ be the space of polynomials on $V$ 
generated by the translated $Q(z-a)$ of $Q$.  We will define a representation
$\kappa$ on $\cal W$ of ${\rm Conf}(V,Q)$ or of a covering of order two of it.

\vfill \eject

\noindent
{\it Case 1}

In case there exists a character $\chi$ of Str$(V,Q)$ such that
$\chi^2=\gamma$, then let $K={\rm Conf}(V,Q)$. Define the cocycle 
$$\mu(g,z)=\chi((Dg(z)^{-1})\quad (g\in K,\ z\in V),$$
and the representation $\kappa$ of $K$ on $\cal W$,
$$(\kappa(g)p)(z)=\mu(g^{-1},z)p(g^{-1}\cdot z).$$
The function $\kappa(g)p$ belongs actually to $\cal W$.
In fact the cocycle $\mu(g,z)$ is a polynomial in $z$ of degree 
$\leq \, {\rm deg}\, Q$ and
$$\eqalign{
(\kappa(\tau_a)p)(z)&=p(z-a) \quad (a\in V),\cr
(\kappa(\ell )p)(z)&=\chi(\ell )p(\ell^{-1}\cdot z) \quad (\ell\in L),\cr
(\kappa(j)p)(z)&=Q(z)p(-z^{-1}).\cr}$$

\medskip

\noindent
{\it Case 2}

Otherwise the group $K$ is defined as the set of pairs
$(g,\mu)$ with $g\in  {\rm Conf}(V,Q)$, and $\mu$ is a rational
function on $V$ such that 
$$\mu(z)^2=\gamma(Dg(z))^{-1}.$$
We consider on $K$ the product  $(g_1,\mu_1)(g_2,\mu_2)=(g_1g_2,\mu_3)$
with
$\mu_3(z)=\mu_1(g_2\cdot z)\mu_2(z)$.
For $\tilde g=(g,\mu) \in K$, define $\mu(\tilde g,z):=\mu(z)$.
Then $\mu(\tilde g,z)$ is a cocycle:
$$\mu(\tilde g_1\tilde g_2,z)=\mu(\tilde g_1,\tilde
g_2\cdot z)\mu(\tilde g_2,z),$$
where $\tilde g\cdot z=g\cdot z$ by definition.

\th Proposition 1.3|
 {\rm (i)} The map 
$$K \rightarrow {\rm Conf}(V,Q),\quad \tilde g=(g,\mu )\mapsto g$$ 
is a surjective group morphism.

{\rm (ii)} For $g\in K$, $\mu(g,z)$ is a polynomial in $z$ of degree
$\leq \,  {\rm deg} \, Q$.
\finth

\proof
It is clearly a group morphism. We will
show that the image contains a set of generators of ${\rm Conf}(V,Q)$.
If $g$ is a translation, then $(g,1)$ and $(g,-1)$ are elements in $K$.
If $g=\ell\in L$, then $Dg(z)=\ell $, and $(\ell ,\alpha), (\ell ,-\alpha)$,
with $\alpha^2=\gamma(\ell )^{-1}$, are elements in $K$.
If $g\cdot z=j(z):=-z^{-1}$, then $Dg(z)^{-1}=P(z)$, where $P(z)$ denotes the
quadratic representation of the Jordan algebra $V$: $P(z)=2T_z^2-T_{z^2}$,  and
$\gamma(P(z))=Q(z)^2$. Then $(j,Q(z)), (j,Q(-z))$ are elements
in $K$.
\hfill \qed

\bigskip

Let $P_{\rm max}$ denote the preimage in $K$ of the maximal parabolic subgroup $L\ltimes N\subset {\rm Conf}(V,Q)$, where $N$ is the subgroup of ranslations. 
For $g\in P_{\rm max}$, $\mu (g,z)$ does not depend on $z$, and $\chi (g)=\mu (g^{-1},z)$ is a character of $P_{\rm max}$. For $g=(\ell ,\alpha )$ ($\ell \in L$), $\chi (g)^2=\gamma (\ell )$.

Observe that the inverse in $K$ of $\sigma=(j,Q(z))$ is
$\sigma^{-1}=(j,Q(-z))$.
If $K$ is connected, then $K$ is a covering of order 2 of
${\rm Conf}(V,Q)$. If not, the identity component $K_0$ of $K$ is
homeomorphic to ${\rm Conf}(V,Q)$.

\vskip 4pt

The representation $\kappa$ of $K$ on $\cal W$ is then given by 
$$\bigl(\kappa(g)p\bigr)(z)=\mu(g^{-1},z)p(g^{-1}\cdot z).$$
In particular
$$\eqalign{
\bigl(\kappa (g)p\bigr)(z)&=\chi (g)p(g^{-1}\cdot z)\quad (g\in P_{\rm max}),\cr
\bigl(\kappa (\sigma )p\bigr)(z)&=Q(-z)p(-z^{-1}).\cr}$$
Hence $p_0\equiv 1$ is a highest weight vector with respect 
and $Q=\kappa (\sigma )p_0$ is a lowest weight vector. 

\bigskip

\noindent
{\it Example 1}

\medskip

If $V={\bboard C}$, $Q(z)=z^n$, then ${\rm Str}(V,Q)={\bboard C}^*$,
$\gamma (\ell )=\ell ^n$, and
${\rm Conf}(V,Q)\simeq PSL(2,{\bboard C})$ is the group of fractional linear transformations
$$z\mapsto g\cdot z={az+b\over cz+d},\ {\rm with}\ g=\pmatrix{a & b \cr c & d \cr}\in SL(2,{\bboard C}).$$
Furthermore
$$Dg(z)={1\over (cz+d)^2},\ \gamma \bigl(Dg(z)^{-1}\bigr)=(cz+d)^{2n},\ \mu (g,z)=(cz+d)^n.$$
Hence, if $n$ is even, then $K=PSL(2,{\bboard C})$,
and, if $n$ is odd, then $K=SL(2,{\bboard C})$.

The space 
$\cal W$ is the space of polynomials of degree $\leq n$ in one variable. The representation
$\kappa $ of $K$ on $\cal W$ is given by
$$\bigl(\kappa (g)p\bigr)(z)=(cz+d)^np\Bigl({az+b\over cz+d}\Bigr),
\ {\rm if}\  g^{-1}=\pmatrix{a & b \cr c & d \cr}.$$ 

\bigskip

\noindent
{\it Example 2}

\medskip

If $V=M(n,{\bboard C})$, $Q(z)=\det z$, then ${\rm Str}(V,Q)=GL(n,{\bboard C})\times GL(n,{\bboard C})$,
acting on $V$ by
$$\ell \cdot z=\ell _1z\ell _2^{-1}\quad \ell =(\ell _1,\ell _2).$$
Then $\gamma (\ell )=\det \ell _1\, \det \ell _2^{-1}$, and $\gamma $ is not the square
of a character of ${\rm Str}(V,Q)$.
Furthermore ${\rm Conf}(V,Q)=PSL(2n,{\bboard C})$
is the group of the rational transformations
$$z\mapsto g\cdot z=(az+b)(cz+d)^{-1},\ {\rm with}\ g=\pmatrix{a & b \cr c & d \cr}\in SL(2n,{\bboard C}),$$
decomposed in $n\times n$-blocs.
To determine the differential of such a transformation, let us write (assuming $c$ to be invertible)
$$g\cdot z=(az+c)(cz+d)^{-1}=
ac^{-1}-(ac^{-1}d-b)(cz+d)^{-1},$$
and we get
$$Dg(z)w=(ac^{-1}d-b)(cz+d)^{-1}cw(cz+d)^{-1}.$$
Notice that $Dg(z)\in {\rm Str}(V,Q)$:
$$Dg(z)w=\ell _1w\ell _2^{-1},\ {\rm with}\  \ell _1=(ac^{-1}d-b)(cz+d)^{-1}c,\ \ell _2=(cz+d).$$
Since $\det (ac^{-1}d-b)\det c=\det g=1$,
$$\gamma \bigl(Dg(z)^{-1} \bigr)=\det (cz+d)^2.$$
It follows that $K=SL(2n,{\bboard C})$, and $\mu (g,z)=\det (cz+d)$.

The space $\cal W$ is a space of polynomials of an $n\times n$ matrix variable,
with degree $\leq n$.
The representation $\kappa $ of $K$ on $\cal W$ is given by
$$\bigl(\kappa (g)p\bigr)(z)=\det (cz+d)p\bigl((az+b)(cz+d)^{-1}\bigr),\
{\rm if}\ g^{-1}=\pmatrix{a & b \cr c & d \cr}.$$

\section 2. The orbit $\Xi $, and the irreducible $K$-invariant Hilbert subspaces 
of ${\cal O}(\Xi )$|
Let $\Xi$  be the
$K$-orbit  of $Q$ in $\cal W$:
$$\Xi=\{\kappa(g)Q \mid g\in K\}.$$
Then $\Xi$ is a conical variety. In fact, if $\xi=\kappa(g)Q$, then, 
for $\lambda\in {\bboard C}^*$, 
$\lambda\xi=\kappa(g\circ h_t)Q$, where $h_t\cdot z=e^{-t}z$ 
($t\in {\bboard C}$)  with $\lambda=e^{2t}$. 

A polynomial $\xi\in \cal W$ can be written
$$\xi(v)=wQ(v)+\ {\rm terms \  of \ degree}\ <4 \quad (w\in {\bboard C}),$$
and $w=w(\xi)$ is a linear form on $\cal W$ which is invariant under the
parabolic subgroup $P_{\rm max}$.
The set $\Xi_0=\{\xi\in \Xi \mid w(\xi)\ne 0\}$
 is open and dense in $\Xi$.
 A polynomial $\xi\in \Xi_0$  can be written
$$\xi(v)=wQ(v-z) \quad
(w \in {\bboard C}^*, z\in V).$$
Hence  we get a coordinate system $(w,z)\in {\bboard C}^*\times V$ for $\Xi_0$.

\th Proposition  2.1| 
In this system, the action
of $K$ is given by 
$$\kappa(g) : (w,z) \mapsto \bigl(\mu(g,z)w,g\cdot z\bigr).$$
\finth

Observe that the orbit $\Xi$ can be seen as a line bundle
over the conformal compactification of $V$.

\bigskip

\proof
Recall that, for $\xi \in \Xi$, 
$$\bigl(\kappa(g)\xi \bigr)(v)=\mu(g^{-1},v)\xi(g^{-1}\cdot v),$$
and, if $\xi(v)=wQ(v-z)$, then
$$=\mu (g^{-1},v)wQ(g^{-1}\cdot v-z)=\mu(g^{-1},v)wQ(g^{-1}\cdot v-g^{-1}g\cdot z).$$
By Lemma 6.6 in [Faraut-Gindikin,1996],
$$\mu(g,z)\mu(g,z')Q(g\cdot z-g'\cdot z')=Q(z-z').$$
Therefore
$$(\kappa(g)\xi)(v)=\mu(g^{-1},g\cdot z)^{-1}wQ(v-g.z)
=\mu(g,z)wQ(v-g\cdot z),$$
by the cocycle property. 
\hfill \qed

\vskip 2pt

The group $K$  acts on the space ${\cal O}(\Xi)$  of holomorphic functions on $\Xi$
by:
$$\bigl(\pi (g)f\bigr)(\xi)=f\bigl(\kappa(g)^{-1}\xi \bigr).$$
If $\xi\in \Xi_0$, i.e.  $\xi(v)=wQ(v-z)$,  and $f\in {\cal O}(\Xi)$, we will write
$f(\xi)=\phi(w,z)$ for the restriction of $f$ to $\Xi_0$.
In the coordinates $(w,z)$, the  representation $\pi$ is given by
$$(\pi(g)\phi)(w,z)=\phi(\mu(g^{-1},z)w,g^{-1}\cdot z).$$
Let ${\cal O}_m(\Xi)$  denote the space of holomorphic functions
$f$ on $\Xi$, homogeneous of degree $m\in {\bboard Z}$: 
$$f(\lambda\xi)=\lambda^mf(\xi)\quad (\lambda\in {\bboard C}^*).$$
The space ${\cal O}_m(\Xi)$  is invariant under the representation
$\pi $. If $f\in {\cal O}_m(\Xi)$, then its 
restriction $\phi$ to $\Xi_0$ can be written $\phi (w,z)=w^m\psi(z)$,
where $\psi $ is a holomorphic function on $V$. We will write 
$\tilde {\cal O}_m(V)$ for the space of the functions $\psi $ corresponding to the
functions $f\in {\cal O}_m(\Xi)$, and denote by $\tilde \pi _m$ the representation of
$K$  on $\tilde {\cal O}_m(V)$  corresponding to the restriction $\pi _m$ of
$\pi$ to ${\cal O}_m(\Xi)$.
The representation $\tilde\pi_m$
is  given by
$$\bigl(\tilde\pi_m(g)\psi \bigr)(z)=\mu(g^{-1},z)^m\psi(g^{-1}\cdot z).$$
Observe that $(\tilde\pi_m(\sigma)1)(z)=Q(-z)^m$.

\th Theorem 2.2|
{\rm (i)} ${\cal O}_m(\Xi)=\{0\}$ for $m<0$.

{\rm (ii)} The space ${\cal O}_m(\Xi)$ is  finite dimensional, and the
representation $\pi_m$ is  irreducible.

{\rm (iii)} The functions $\psi$ in $\tilde{\cal O}_m(V)$  are polynomials.
\finth

\proof
(i)  Assume ${\cal O}_m(\Xi )\not=\{0\}$. Let $f\in {\cal O}_m(\Xi)$, $f\not\equiv 0$, and $\phi(w,z)=\psi(z)w^m$ its
restriction to $\Xi_0$. Then $\psi $ is holomorphic on $V$, and
$$\bigl(\tilde\pi_m(\sigma)\psi \bigr)(z)=Q(-z)^m\psi(-z^{-1}),$$
is holomorphic as well. We may assume $\psi (e)\not=0$. The function $h(\zeta)=\psi (\zeta e) \quad (\zeta \in {\bboard C})$
is holomorphic on $\bboard C$,
$$h(\zeta)=\sum _{k=0}^{\infty} a_k\zeta^k,$$
together with the function 
$$Q(\zeta e)^m\psi(-{1\over \zeta}e)=\zeta^{mN}
h(-{1\over \zeta})=\zeta^{mN} \sum _{k=0}^{\infty}a_k(-{1\over \zeta})^k
\quad (N={\rm deg}\, Q).$$
It follows that $m\geq 0$, and that $a_k=0$ for $k>mN$.

(ii)  The subspace 
$$\{f\in {\cal O}_m(\Xi)\mid \forall a\in V, \pi(\tau_a)f=f\}$$
reduces to the functions $Cw^m$, hence is one dimensional. By the theorem of the
highest weight [Goodman,2008], it follows that ${\cal O}_m(\Xi)$ is finite
dimensional and irreducible.

(iii)  Furthermore it follows that the functions in ${\cal O}_m(\Xi)$ are of the
form $w^m\psi(z)$,
where $\psi$ is a polynomial on $V$ of degree $\leq m\cdot\deg Q$.
\hfill \qed

\vskip 6pt
We fix a  Euclidean real form $V_{\bboard R}$ of the complex Jordan algebra $V$, denote by $z\mapsto
\bar z$ the conjugation of $V$ with respect to $V_{\bboard R}$, and then 
consider  the involution $g\mapsto \bar g$ of ${\rm Conf}(V,Q)$ given  by: 
$\bar g\cdot z=\overline{g\cdot \bar z}$.  
For $(g,\mu)\in K$ define
$$\overline{(g,\mu)}=(\bar g,\bar\mu),\  {\rm where}\ 
\bar\mu(z)=\overline{\mu(\bar z)}.$$

The involution
$\alpha$ defined by $\alpha (g)=\sigma\circ\bar g\circ\sigma^{-1}$ is a Cartan
involution of $K$ (see Proposition 1.1. in  [Pevzner,2002]),  and  
$$K_{\bboard R}:=\{g\in K\mid \alpha(g)=g\}$$
is a compact real form of $K$.
 
\bigskip

\noindent
{\it Example 1.}

\medskip

If $V={\bboard C}$, $Q(z)=z^n$. Then $V_{\bboard R}={\bboard R}$, and
$z\mapsto \bar z$ is the usual conjugation.  We saw that
$K=PSL(2,{\bboard C})$ if $n$ is even, and $SL(2,{\bboard C})$ if $n$ is odd. For $g\in SL(2,{\bboard C})$,
$$g=\pmatrix{
a & b \cr
c & d \cr},$$
we get 
$$\alpha (g)
=\pmatrix{
0 & 1 \cr
-1 & 0 \cr}
\pmatrix{ 
\bar a&\bar b \cr
\bar c & \bar d \cr}
\pmatrix{
0 & -1 \cr
1 & 0 \cr}
=\pmatrix{
\bar d &-\bar c \cr
-\bar b & \bar a \cr}.$$
Hence $K_{\bboard R}=PSU(2)$ if $n$ is even, and $K_{\bboard R}=SU(2)$ if $n$ is odd.

\bigskip

\noindent

{\it Example 2.} 

\medskip

If $V=M(n,{\bboard C})$, $Q(z)=\det z$,
then $V_{\bboard R}=Herm(n,{\bboard C})$ and the conjugation is $z\mapsto z^*$. 
We saw that $K=SL(2n,{\bboard C})$.
For $g\in SL(2n,{\bboard  C})$,
$$g=\pmatrix{
a & b \cr
c & d \cr} ,$$
we get 
$$\alpha(g)=
\pmatrix{
0 & I \cr
-I & 0 \cr}
\pmatrix{
a^* & b^* \cr
c^* & d^* \cr}
\pmatrix{
0 & -I \cr
I &0 \cr}
=\pmatrix{
d^* & -c^* \cr
-b^* & a^*\cr}.$$
Hence $K_{\bboard R}=SU(2n)$.

\bigskip

We will define on ${\cal O}_m(\Xi )$ a $K_{\bboard R}$-invariant inner product.
Define the subgroup $K_0$ of $K$ as $K_0=L$ in Case 1, and the preimage
of $L$ in Case 2, relatively to the covering map $K\to {\rm Conf}(V,Q)$,
and also $(K_0)_{\bboard R}=K_0 \cap K_{\bboard R}$. 
The coset space $M=K_{\bboard R}/(K_0)_{\bboard R}$, 
is a compact  Hermitian space and is the
conformal compactification of $V$. 
There is on $M$ a $K_{\bboard R}$-invariant probability measure, 
for which $M\setminus V$ has measure 0. Its restriction $m_0$ to $V$ is a probability measure with a density which can be computed by using the decomposition of $V$ into simple
Jordan algebras.

Let $H(z,z')$ be the polynomial on $V\times V$, holomorphic in $z$,
anti-holomorphic in $z'$ such that 
$$H(x,x)=Q(e+x^2) \quad (x\in V_{\bboard R}).$$
Put $H(z)=H(z,z)$. If $z$ is invertible, then
$H(z)=Q(\bar z)Q(\bar z^{-1}+z)$.

\th Proposition 2.3|
For $g\in K_{\bboard R}$,
$$H(g\cdot z_1,g\cdot z_2)\mu(g,z_1)\overline{\mu(g,z_2)}=H(z_1,z_2),$$
and 
$$H(g\cdot z)\vert\mu(g,z)\vert^2=H(z).$$
\finth

\proof
Recall that an element $g\in K_{\bboard R}$ satisfies $\sigma\circ
\bar g\circ\sigma^{-1}=g$, or $\sigma\circ\bar g=g\circ \sigma$. Recall also the cocycle
property: for $g_1,g_2 \in K$,
$$\mu(g_1g_2,z)=\mu(g_1,g_2\cdot z)\mu(g_2,z).$$
Since $\mu(\sigma,z)=Q(z)$, it follows that, for $g\in K_{\bboard R}$,
$$\mu(g,\sigma\cdot z)Q(z)=Q(\bar g\cdot z)\mu(\bar g,z). \eqno (1)$$
By Lemma 6.6 in [Faraut-Gindikin,1996], for $g\in K$,
$$Q(g\cdot z_1-g\cdot z_2)\mu(g,z_1)\mu(g,z_2)=Q(z_1-z_2). \eqno (2)$$
For $g\in K_{\bboard R}$,
$$\eqalign{
H(g\cdot z_1,g\cdot z_2)&=Q(\bar g\cdot z_2)Q(g\cdot z_1-\sigma\bar g\cdot \bar z_2)\cr
&=Q(\bar g \cdot \bar z_2)Q(g\cdot z_1-g\sigma\bar z_2),\cr}$$
and, by (2),
$$=Q(\bar g\cdot\bar z_2)\mu(g,z_1)^{-1}\mu(g,\sigma\cdot\bar z_2)^{-1}
Q(z_1-\sigma\cdot\bar z_2).$$
Finally, by (1),
$$=\mu(g,z_1)^{-1}\mu(\bar g,\bar z_2)^{-1}H(z_1,z_2). \eqno \qed $$

\bigskip

We define the norm of a function $\psi\in \tilde{\cal O}_m(V)$ by
$$\Vert \psi\Vert_m^2={1\over a_m}
\int_V\vert \psi(z)\vert^2H(z)^{-m}m_0(dz),$$
with 
$$a_m=\int_VH(z)^{-m}m_0(dz).$$

\th Proposition 2.4|
{ \rm (i)}  This norm is $K_{\bboard R}$-invariant. 
Hence,  $\tilde{\cal O}_m(V)$ is a Hilbert subspace of
${\cal O}(V)$. 

{\rm (ii)} The reproducing kernel of $\tilde{\cal O}_m(V)$ is given by
$$\tilde{\cal K}_m(z,z')= H(z,z')^m.$$
\finth

\proof
(i)  From Proposition 2.3 it follows that, for $g \in K_{\bboard R}$,
$$\eqalign{
\Vert\tilde\pi_m(g^{-1})\psi\Vert_m^2&
={1\over a_m}\int_V \vert\mu(g,z)\vert^{2m}
\vert\psi (g^{-1}\cdot z)\vert^2H(z)^{-m}m_0(dz) \cr
&={1\over a_m}\int_V\vert\psi(g^{-1}\cdot z)\vert^2H(g^{-1}\cdot z)^{-m}m_0(dz) \cr
&={1\over a_m}\int_V\vert\psi(z)\vert^2H(z)^{-m}m_0(dz)
=\Vert\psi\Vert_m^2. \cr}$$

(ii) There is a unique function $\psi_0 \in  \tilde{\cal O}_m(V)$ such
that, for  $\psi \in \tilde{\cal O}_m(V)$,
$$(\psi\mid\psi_0)=\psi (0).$$
The function $\psi_0$ is $K_0$-invariant, therefore constant:
$\psi _0(z)=C$. Taking $\psi=\psi_0$, one gets  $C^2=C$, hence $C=1$. It
means that, if $\tilde{\cal K}_m(z,z')$ denotes the reproducing kernel of
$\tilde{\cal O}_m(V)$,
$$\tilde{\cal K}_m(z,0)=\tilde{\cal K}_m(0,z')=1.$$
Since $\tilde{\cal K}_m(z,z')$ and $H(z,z')$ satisfy the following
invariance properties: for $g\in K_{\bboard R}$,
$$\eqalign{
\tilde{\cal K}_m(g\cdot z,g\cdot z')\mu(g,z)^m\overline{\mu(g,z')}^m
&=\tilde{\cal K}_m(z,z'), \cr
H(g\cdot z,g\cdot z')\mu(g,z)\overline{\mu(g,z')}&=H(z,z'),\cr}$$
it follows that
$$\tilde{\cal K}_m(z,z')=H(z,z')^m . \eqno \qed $$
\bigskip

Since ${\cal O}_m(\Xi)$ is isomorphic to $\tilde{\cal O}_m(V)$, 
the space ${\cal O}_m(\Xi)$ becomes an invariant Hilbert subspace of ${\cal O}(\Xi)$, with
reproducing kernel 
$${\cal K}_m(\xi,\xi')=\Phi (\xi,\xi')^m,$$
where 
$$\Phi (\xi,\xi')=H(z,z')w\overline{w'} \qquad
(\xi=(w,z), \xi'=(w',z')).$$

\th Theorem 2.5|
The group $K_{\bboard R}$  acts
multiplicity free on ${\cal O}(\Xi)$.  
The irreducible $K_{\bboard R}$-invariant  subspaces of ${\cal O}(\Xi)$ 
are the spaces ${\cal O}_m(\Xi)$ ($m\in {\bboard N}$).  
If ${\cal H} \subset {\cal O}(\Xi)$ is a 
$K_{\bboard R}$-invariant Hilbert subspace, the reproducing kernel
of $\cal H$ can be written 
$${\cal K}(\xi,\xi')=\sum _{m=0}^{\infty}c_m\Phi (\xi,\xi')^m,$$
with $c_m\geq 0$, such that the series 
$\sum _{m=0}^{\infty} c_m\Phi(\xi,\xi')^m$ converges uniformly on compact subsets in $\Xi$.
\finth

This multiplicity free property means that $K_{\bboard R}$ acts multiplicity
free on every $K_{\bboard R}$-invariant Hilbert space ${\cal H}\subset {\cal O}(\Xi )$.

\bigskip

\proof
The representation $\pi$ of $K_{\bboard  R}$ on ${\cal O}(\Xi)$
commutes with the ${\bboard C}^*$-action by dilations and the spaces ${\cal O}_m(\Xi)$ are
irreducible, and mutually inequivalent. It follows that $K_{\bboard R}$ acts multiplicity free. 
\hfill \qed

\bigskip

In case of a weighted Bergman space there is an integral formula for the numbers $c_m$.
For a positive function $p(\xi)$ on $\Xi$, 
consider the subspace ${\cal H} \subset {\cal O}(\Xi)$ 
of functions $\phi$ such that 
$$\Vert\phi\Vert^2=\int_{{\bboard C}\times V}\vert\phi(w,z)\vert^2p(w,z)m(dw)m_0(dz) <\infty,$$
where $m(dw)$ denotes the Lebesgue measure on $\bboard C$.

\th Theorem 2.6|
Let $F$ be a positive function on
$[0,\infty[$, and define 
$$p(w,z)=F(H(z)\vert w\vert^2)H(z).$$

{\rm (i)} Then $\cal H$ is $K_{\bboard R}$-invariant.

{\rm (ii)} If
$$\phi(w,z)=\sum _{m=0}^{\infty}w^m\psi_m(z),$$
then
$$\Vert
\phi\Vert^2=\sum _{m=0}^{\infty}{1\over c_m}\Vert\psi_m\Vert_m^2,$$
with
$${1\over c_m}=\pi a_m \int_0^{\infty}F(u)u^mdu.$$

{\rm (iii)} The reproducing kernel of $\cal H$ is given by 
$${\cal K}(\xi,\xi')=\sum _{m=0}^{\infty}c_m\Phi(\xi,\xi')^m.$$
\finth

\proof
a) Observe first that the function defined on $\Xi$ by 
$$(w,z) \mapsto \vert w\vert^2H(z),$$
is $K_{\bboard R}$-invariant. In fact, for $g\in K$,
$$\kappa(g) : (w,g) \mapsto \bigl(\mu(g,z)w,g\cdot z\bigr),$$ 
and, by Propositiion 2.3, for $g\in K_{\bboard R}$,
$$\vert\mu(g,z)\vert^2H(g\cdot z)=H(z).$$
Furthermore the measure $h(z)m(dw)m_0(dz)$ is also invariant under $K_{\bboard R}$.
In fact, under the transformation $z=g\cdot z' , w=\mu(g,z')w'$ ($g\in
K_{\bboard R}$), we get
$$\eqalign{
H(z)m(dw)m_0(dz) 
&=H(g\cdot z')\vert\mu(g,z')\vert^2m(dw')m_0(dz') \cr
&=H(z')m(dw')m_0(dz').\cr}$$

b) Assume that $p(w,z)=F\bigl(H(z)\vert w\vert^2\bigr)H(z)$. Then
$$\Vert\pi(g)\phi\Vert^2=\int_{{\bboard C}\times V}
\vert\phi(\mu(g^{-1},z)w,g^{-1}\cdot z)\vert^2
F\bigl(H(z)\vert w\vert^2\bigr)H(z)m(dw)m_0(dz).$$
We put
$$g^{-1}\cdot z=z'\quad , \quad \mu(g^{-1},z)w=w'.$$
By the invariance of the measure $H(z)m(dw)m_0(dz)$, we obtain
$$\eqalign{
&\Vert\pi(g)\phi\Vert^2 =\cr
&\int_{{\bboard C}\times V}
\vert\phi(w',z')\vert^2F\bigl(H(g\cdot z')\vert\mu(g^{-1},g\cdot z')\vert^{-2}
\vert w'\vert^2\bigr)H(z')m(dw')m_0(dz').\cr}$$
Furthermore 
$$H(g\cdot z')\vert\mu(g^{-1},g\cdot z')\vert^{-2}
=H(g\cdot z')\vert\mu(g,z')\vert^2=H(z'),$$ 
and, finally,  
$\Vert\pi(g)\phi\Vert=\Vert\phi\Vert$.

c) If $\phi(w,z)=w^m\psi(z)$, then
$$\Vert\phi\Vert^2=\int_{{\bboard C}\times V}\vert
w\vert^{2m}\vert\psi(z)\vert^2F\bigl(H(z)\vert w\vert^2\bigr)H(z)m(dw)m_0(dz').$$
We put $w'=\sqrt{H(z)}w$, then
$$\eqalignno{
\Vert\phi\Vert^2&=\int_{{\bboard C}\times V}H(z)^{-m}
\vert w'\vert^{2m}\vert\psi(z)\vert^2F(\vert w'\vert^2)m(dw')m_0(dz') \cr
&=a_m\Vert\psi\Vert_m^2\int_{\bboard C}F(\vert w'\vert^2)\vert w'\vert^{2m}m(dw') \cr
&=a_m\Vert\psi\Vert_m^2\pi\int_0^{\infty}F(u)u^mdu. & \qed \cr}$$

\section 3. Decomposition into simple Jordan algebras|
Let us decompose the semi-simple  Jordan algebra $V$ into simple ideals:
$$V=\bigoplus _{i=1}^s V_i.$$
Denote by $n_i$ and $r_i$ the dimension and the rank of the simple Jordan algebra $V_i$, and $\Delta _i$ the determinant polynomial. Then
$$Q(z)=\prod _{i=1}^s \Delta _i(z_i)^{k_i}.$$
Let $H_i(z,z')$ be the polynomial on $V_i\times V_i$, holomorphic in $z$, antiholomorphic in $z'$, such that
$$H_i(z,z')=\Delta _i(e_i+x^2) \quad \bigl(x\in (V_i)_{\bboard R}\bigr),$$
and put $H_i(z)=H_i(z,z)$. The measure $m_0$ has a density with respect to the Lebesgue measure $m$ on $V$:
$$m_0(dz)={1\over C_0}H_0(z) m(dz),$$
with
$$\eqalign{
H_0(z)&=\prod _{i=1}^s H_i(z_i)^{-2{nÐi\over r_i}}, \cr
C_0&=\int _V H_0(z)m(dz). \cr}$$
The Lebesgue measure $m$ will be chosen such that $C_0=1$.

\th Proposition 3.1|
{\rm (i)} The polynomial $Q$ satisfies the following Bernstein identity
$$Q\Bigl({\partial \over \partial z}\Bigr) Q(z)^{\alpha } =B(\alpha )Q(z)^{\alpha -1} 
\quad (z\in {\bboard C}),$$
where the Bernstein polynomial $B$ is given by
$$B(\alpha )=\prod _{i=1}^s b_i(k_i\alpha )b_i(k_i\alpha -1)\ldots b_i(k_i\alpha -k_i+1),$$
and $b_i$ is the Bernstein polynomial relative to the determinant polynomial $\Delta _i$.

{\rm (ii)} Furthermore
$$Q\Bigl({\partial \over \partial z}\Bigr) H(z)^{\alpha }
=B(\alpha )\overline{Q(z)}H(z)^{\alpha -1}.$$
\finth

\proof
(i) The Bernstein identity for $Q$ follows from Proposition VII.1.4 in [Faraut-Kor\'anyi,1994].

(ii) For $z$ invertible
$$H(z)=Q(\bar z)Q(\bar z^{-1}+z),$$
and then, by (i),
$$\eqalignno{
Q\Bigl({\partial \over \partial z}\Bigr)H(z)^{\alpha }
&=Q(\bar z)^{\alpha }B(\alpha )Q(\bar z^{-1}+z)^{\alpha -1}\cr
&=Q(\bar z)B(\alpha )H(z)^{\alpha -1}.& \qed \cr}$$

\bigskip

\noindent
{\it Example 1}

\medskip

If $V={\bboard C}$, $Q(z)=z^n$, then
$$\Bigl({d\over dz}\Bigr)^n z^{n\alpha }=B(\alpha )z^{n(\alpha -1)},$$
with
$$B(\alpha )=n\alpha (n\alpha -1)\ldots (n\alpha -n+1).$$

\vfill \eject

\noindent
{\it Example 2}

\medskip

If $V=M(n,{\bboard C})$, $Q(z)=\det z $, then
$$\det \Bigl({\partial \over \partial z}\Bigr)(\det z)^{\alpha }=B(\alpha )(\det z)^{\alpha -1},$$
with
$$B(\alpha )=\alpha (\alpha+1)\ldots (\alpha +n-1).$$

\bigskip

Recall that we have introduced the numbers
$$a_m=\int _V H(z)^{-m}m_0(dz).$$

\tha Proposition 3.2|
$$a_m=\prod _{i=1}^s {\Gamma _{\Omega _i}(2{n_i\over r_i})
\over \Gamma _{\Omega _i}({n_i\over r_i})}
\prod _{i=1}^s {\Gamma _{\Omega _i}(mk_i+{n_i\over r_i})
\over \Gamma _{\Omega _i}(mk_i+2{n_i\over r_i})},$$
where $\Gamma _{\Omega _ i}$ is the Gindikin gamma function of the symmetric cone $\Omega _i$ in the Euclidean Jordan algebra $(V_i)_{\bboard R}$.
\fintha

\proof
If the Jordan algebra  $V$ is simple and $Q=\Delta $, the determinant polynomial,
by Proposition X.3.4 in [Faraut-Kor\'anyi,1994],
$$\eqalign{
a_m
&=\int _V H(z)^{-m}m_0(dz)={1\over C_0}\int _V H(z)^{-m-2{n\over r}}m(dz) \cr
&=C\int _{\Omega } \Delta (e+x)^{-m-2{n\over r}}m(dx).\cr}$$
By Exercice 4 of Chapter VII in [Faraut-Kor\'anyi,1994] we obtain
$$a_m=C'{\Gamma _{\Omega }(m+{n\over r})\over \Gamma _{\Omega } (m+2{n\over r})}.$$
In the general case
$$a_m={1\over C_0}\prod _{i=1}^s \int _{V_i}H_i(z_i)^{-mk_i-2{n_i\over r_i}}m_i(dz_i),$$
and the formula of the proposition follows.
\hfill \qed

\vfill \eject

\section 4.  Generalized Kantor--Koecher--Tits construction|
From now on,  $Q$ is assumed to be of degree $4$.
The group of dilations of $V$ : $h_t\cdot z=e^{-t}z$ ($t \in {\bboard C}$) is a
one parameter subgroup of $L$, and $\chi(h_t)=e^{-2t}$.  Put
$h_t=\exp (tH)$. Then $\ad (H)$ defines a grading of the Lie algebra
$\goth k$ of $K$:
$${\goth k}={\goth k}_{-1}+{\goth k}_0+{\goth k}_1,$$
with ${\goth k}_j=\{X\in {\goth k}\mid \ad (H)X=jX\}$, ($j=-1,0,1$). Notice that 
$${\goth k}_{-1}=Lie(N)\simeq V,\quad {\goth k}_0=Lie(L),\quad
\Ad (\sigma): {\goth k}_j\rightarrow {\goth k}_{-j},$$
and also that  $H$ belongs to the centre ${\goth z}({\goth k}_0)$ of ${\goth k}_0$.
The element $H$ defines also  a grading of ${\goth p}:={\cal W}$:
$${\goth p}={\goth p}_{-2}+{\goth p}_{-1}+{\goth p}_0+{\goth p}_1+{\goth p}_2,$$
where   
$${\goth p}_j=\{p\in {\goth p} \mid d\kappa(H)p=jp\}$$
is the set of polynomials in $\goth p$,  homogeneous  of degree 
$j+2$. The subspaces ${\goth p}_j$ are invariant under $K_0$. 
Furthermore 
$\kappa(\sigma) : {\goth p}_j \rightarrow
{\goth p}_{-j}$, and
$${\goth p}_{-2}={\bboard C}, \quad {\goth p}_2={\bboard C}\, Q, 
\quad {\goth p}_{-1} \simeq V,\quad {\goth p}_1\simeq V.$$
Let ${\goth g}={\goth k}\oplus {\goth p}$. Put $E=Q$, $F=1$.

\th Theorem 4.1|
There exists on  $\goth g$ a 
unique Lie algebra structure such that: 
$$\eqalign{
(i) \quad [X,X']&=[X,X']_{\goth k} \quad (X,X' \in {\goth k}), \cr
(ii) \quad [X,p]&=d\kappa(X)p \quad (X\in {\goth k}, p\in {\goth p}), \cr
(iii) \quad [E,F]&=H. \cr}$$
\finth

\proof
Observe that $(E,F,H)$ is an ${\goth sl}_2$-triple, and 
that $H$ defines a grading of 
$${\goth g}={\goth g}_{-2}+{\goth g}_{-1}+{\goth g}_0+{\goth g}_1+ {\goth g}_2,$$
with 
$${\goth g}_{-2}={\goth p}_{-2}, \quad
{\goth g}_{-1}={\goth k}_{-1}+{\goth p}_{-1},\quad
{\goth g}_0={\goth k}_0+{\goth p}_0,\quad
{\goth g}_1={\goth k}_1+{\goth p}_1, \quad
{\goth g}_2={\goth p}_2.$$
It is possible  to give a direct proof of Theorem 4.1 (see  Theorem 3.1. in 
[Achab,2011]). It is also possible to see this statement as a special
case of constructions of Lie algebras by Allison and Faulkner [1984]. We
describe below this construction in our case. 

\vfill \eject

a) {\it Cayley-Dickson process.}   

Let $x\mapsto x^*$
denote the symmetry with respect to the one dimensional subspace ${\bboard C} e$:
$$x^*=\demi \langle x,e\rangle -x.$$
Observe that
$$\langle x,e\rangle =\tau (T_x)=D_x\log Q(e),\quad \langle e,e\rangle =4.$$
On the vector space $W=V\oplus V$, one defines an algebra structure: if $z_1=(x_1,y_1),
z_2=(x_2,y_2)$, then $z_1z_2=z=(x,y)$ with 
$$x=x_1x_2-(y_1y_2^*)^* ,\quad  y=x_1^*y_2+(y_1^*x_2^*)^*,$$
and an involution 
$$\bar z=\overline{(x,y)}=(x,-y^*).$$
This involution is an antiautomorphism: 
$\overline{z_1z_2}=\bar z_2\bar z_1$. 
For $a,b\in W$, one introduces the endomorphisms $V_{a,b}$ and
$T_a$ given by 
$$\eqalign{
V_{a,b}z&=\{a,b,z\}:=(a\bar b)z+(z\bar b)a-(z\bar a)b,\cr
T_az&=V_{a,e}z=az+z(a-\bar a).\cr}$$
By Theorem 6.6 in [Allison-Faulkner, 1984] the algebra $W$ is
structurable. This means that, for $a,b,c,d\in W$, 
$$[V_{a,b},V_{c,d}]=V_{V_{a,b}c,d}-V_{c,V_{b,a}d}. \eqno (*)$$
Moreover the structurable algebra $W$ is simple. By $(*)$, the vector
space spanned by the endomorphisms $V_{a,b}$ ($a,b\in W$) is a Lie
algebra denoted by $Instrl(W)$. This algebra is the Lie algebra
${\goth g}_0$ in the grading, and its subalgebra ${\goth k}_0$ is the
structure algebra of the Jordan algebra $V$.
The space $S$ of skew-Hermitian elements in $W$, 
$S=\{z\in W \mid \bar z=-z\}$, 
has dimension one. Its elements are proportionnal to $s_0=(0,e)$.
The subspace $\{(x,0)\mid x\in V\}$ of $W$ is identified to $V$, and
any element $z=(x,y) \in W$ can be written $z=x+s_0y$.

\medskip

b) {\it Generalized Kantor-Koecher-Tits construction.}

One defines a bracket on the vector space
$${\cal K}(W)=\tilde S\oplus\tilde W\oplus Instrl(W)\oplus W\oplus S,$$
where $\tilde S$ is a second copy of $S$, and $\tilde W$ of $W$.
This construction is described in
[Allison,1979], and, by Corollary 6 in that paper, ${\cal K}(W)$ is a
simple Lie algebra. On the subspace ${\cal K}(V)=\tilde V\oplus {\goth str}(V)\oplus V$,
this construction  agrees with the classical Kantor-Koecher-Tits
construction, which  produces   the Lie algebra 
${\goth k}={\goth k}_{-1}\oplus {\goth k}_0\oplus {\goth k}_1$.
This algebra ${\cal K}(W)$ satisfies property (i): the
restriction of the bracket of ${\cal K}(W)$ to ${\cal K}(V)$ coincides
to the one of ${\cal K}(V)$. It satisfies (iii) as well: $[s_0,\tilde
s_0]=I$, the identity of $End(W)$. It remains to check property
(ii). This can be seen as a consequence of the theorem of the
highest weight for irreducible  finite dimensional representations
of reductive Lie algebras. In fact, the representation $d\kappa$ of
$\goth k$ on $\goth p$ is irreducible with highest weight vector $Q$,
with respect to any Borel subalgebra 
${\goth b}\subset {\goth k}_0+{\goth k}_1$ : 

- If $X\in {\goth k}_1$, then d$\kappa(X)Q=0$.

- If $X\in {\goth k}_0$, such that $d\gamma(X)=0$, then
$d\kappa(X)Q=0$, and $d\kappa(H)Q=2Q$. 

On the other hand, for the bracket of ${\cal K}(W)$,

- If $u\in V, [u,s_0]=0$.

- If $X\in {\goth str}(V)$, such that $\tr (X)=0$, then $[X,s_0]=0$ and
$[H,s_0]=2s_0$.

It follows that the adjoint representation of 
${\cal K}(V)=\tilde V\oplus {\goth str}(V)\oplus V$ on 
$$\tilde S\oplus \tilde s_0\tilde V\oplus T_W\oplus s_0V\oplus S,$$
where $T_W=\{T_w=V_{w,e} \mid w\in W\}$,
agrees with the representation $d\kappa$ of ${\goth k}$ on ${\goth p}$.
 In the present case, $T_w=L(w)+\demi \langle v,e\rangle Id$, if
$w=u+s_0v$ ($u,v\in V$).

On the vector space
$${\goth g}={\goth g}_{-2}\oplus {\goth g}_{-1}
\oplus {\goth g}_0\oplus {\goth g}_1\oplus {\goth g}_2,$$
with 
$${\goth g}_1=W ,\quad {\goth g}_{-1}=W ,\quad {\goth g}_2={\bboard C}\, E ,
\quad {\goth g}_{-2}={\bboard C}\, F,\quad {\goth g}_0= Instrl(W),$$
one defines a bracket satisfying the following properties:

(1) ${\goth g}_1+{\goth g}_2$ is a Heisenberg Lie algebra:
$${\goth g}_1 \times {\goth g}_1 \rightarrow {\goth g}_2,\quad
(w_1,w_2)\mapsto w_1\bar w_2-w_2\bar w_1=\psi(w_1,w_2)s_0.$$
The bilinear form $\psi$ is skew symmetric, and $[w_1,w_2]=\psi(w_1,w_2)E$.

(2) ${\goth g}_1\times {\goth g}_{-1}\rightarrow {\goth g}_0, \quad (w,\tilde
w)\mapsto V_{w,\tilde w}$.

(3) ${\goth g}_2 \times {\goth g}_{-1} \rightarrow {\goth g}_1, \quad 
(\lambda E,\tilde w)\mapsto \lambda \tilde w$.
\hfill \qed

\vfill \eject

We introduce now a real form ${\goth g}_{\bboard R}$ of $\goth g$ which will be considered in the sequel.
In Section 2 we have considered the involution $\alpha $ of $K$ given by
$$\alpha (g)=\sigma \circ \bar g \circ \sigma ^{-1}\quad (g\in K),$$
and the compact real form $K_{\bboard R}$ of $K$:
$$K_{\bboard R}=\{g\in K\mid \alpha (g)=g\}.$$
Recall that ${\goth p}$ has been defined as a space of polynomial functions on $V$. For $p\in {\goth p}$, define
$$\bar p=\overline{p(\bar z)},$$
and consider the antilinear involution $\beta $ of $\goth p$ given by
$$\beta (p)=\kappa (\sigma )\bar p.$$
Observe that $\beta (E)=F$.
The involution $\beta $ is related to the involution $\alpha $ of $K$ by the relation
$$\kappa \bigl(\alpha (g)\bigr) \circ \beta =\beta \circ \kappa (g) \quad (g\in K).$$
Hence, for $g\in K_{\bboard R}$, $\kappa (g)\circ \beta =\beta \circ \kappa (g)$. Define
$${\goth p}_{\bboard R}=\{p\in {\goth p}\mid \beta (p)=p\}.$$
The real subspace ${\goth p}_{\bboard R}$ is invariant under $K_{\bboard R}$, and irreducible for that action. The space $\goth p$, as a real vector space, 
decomposes under $K_{\bboard R}$ into two irreducible subspaces
$${\goth p}={\goth p}_{\bboard R}\oplus i{\goth p}_{\bboard R}.$$
One checks that $E+F\in {\goth p}_{\bboard R}$ (and hence $i(E-F)$ as well).

Let $\goth u$ be a compact real form of $\goth g$ such that 
${\goth k}\cap {\goth u}={\goth k}_{\bboard R}$, the Lie algebra of $K_{\bboard R}$. Then $\goth p$ decomposes as
$${\goth p}={\goth p}\cap (i{\goth u})\oplus {\goth p}\cap {\goth u}$$
into two irreducible $K_{\bboard R}$-invariant real subspaces. Looking at the subalgebra ${\goth g}^0$ isomorphic to ${\goth sl}(2,{\bboard C})$ generated by the triple $(E,F,H)$, one sees that $E+F\in {\goth p}\cap (i{\goth u})$. Therefore
${\goth p}_{\bboard R}={\goth p}\cap (i{\goth u})$, and
$${\goth g}_{\bboard R}={\goth k}_{\bboard R}\oplus {\goth p}_{\bboard R}$$
is a Lie algebra, real form of $\goth g$, and the above decomposition is a Cartan decomposition of ${\goth g}_{\bboard R}$.

\vskip 1cm

For the table of next page we have used the notation: 
$$\varphi _n(z)=z_1^2+\cdots +z_n^2,\quad (z\in {\bboard C}^n).$$
In case of an exceptional Lie algebra $\goth g$, the real form ${\goth g}_{\bboard R}$ has been identified by computing the Cartan signature.

\vfill \eject

\hsize=16cm

\def\tvi{\vrule height 12pt depth 1pt width 0pt}

\def\traithorizontal{\noalign{\hrule}}

$$\vbox{\offinterlineskip \halign{
\tvi # & # & # & # & # \cr
$V$ & $Q$ & $\goth k$ & $\goth g$ & ${\goth g}_{\bboard R}$ \cr
\tvi \cr
\traithorizontal 
${\bboard C}^n$ & $\varphi _n(z)^2$ & ${\goth so}(n+2,{\bboard C})$ \hfill 
& ${\goth sl}(n+2,{\bboard C})$ 
& ${\goth sl}(n+2,{\bboard R})$ \cr
\tvi \cr
\traithorizontal
${\bboard C}^p\oplus {\bboard C}^q$ \hfill & $\varphi _p(z)\varphi _q(z')$
& ${\goth so}(p+2,{\bboard C})\oplus {\goth so}(q+2,{\bboard C})$
& ${\goth so}(p+q+4,{\bboard C})$ & ${\goth so}(p+2,q+2)$ \cr
\tvi \cr
\traithorizontal
$Sym(4,{\bboard C})$ & $\det z$ & ${\goth sp}(8,{\bboard C})$ & ${\goth e}_6$
& ${\goth e}_{6(6)}$ \cr
\tvi \cr
$M(4,{\bboard C})$ & $\det z$ & ${\goth sl}(8,{\bboard C})$ 
& ${\goth e}_7$ & ${\goth e}_{7(7)}$ \cr
\tvi \cr
$Skew(8,{\bboard C})$ & ${\rm Pfaff}(z)$ & ${\goth so}(16,{\bboard C})$ & ${\goth e}_8$
& ${\goth e}_{8(8)}$ \cr
\tvi \cr
\traithorizontal
$Sym(3,{\bboard C})\oplus {\bboard C}$ & $\det z \cdot z'$ & 
${\goth sp}(6,{\bboard C})\oplus {\goth sl}(2,{\bboard C})$ \hfill 
& ${\goth f}_4$ & ${\goth f}_{4(4)}$ \cr
\tvi \cr
$M(3,{\bboard C})\oplus {\bboard C}$ \hfill & $\det z\cdot z' $ & 
${\goth sl}(6,{\bboard C})\oplus {\goth sl}(2,{\bboard C})$ \hfill &
${\goth e}_6$ & ${\goth e}_{6(2)}$ \cr
\tvi \cr
$Skew(6,{\bboard C})\oplus {\bboard C}$ & ${\rm Pfaff}( z) \cdot z' $ 
& ${\goth so}(12,{\bboard C})\oplus {\goth sl}(2,{\bboard C})$ \hfill & 
${\goth e}_7$ & ${\goth e}_{7(-5)}$ \cr
\tvi \cr
$Herm(3,{\bboard O})_{\bboard C}\oplus {\bboard C}$ & $\det z\cdot z'$ &
${\goth e}_7\oplus {\goth sl}(2,{\bboard C})$ \hfill & ${\goth e}_8$ & ${\goth e}_{8(-24)}$ \cr
\tvi \cr
\traithorizontal
${\bboard C}\oplus {\bboard C}$ \hfill & $z^3\cdot z'$ & 
${\goth sl}(2,{\bboard C})\oplus {\goth sl}(2,{\bboard C})$ \hfill & ${\goth g}_2$ & ${\goth g}_{2(2)}$ \cr
\tvi \cr
\traithorizontal
}}$$

\vfill \eject

\hsize=11,25cm

\section 5. Representation of the generalized
Kantor-Koecher-Tits Lie algebra| 
Following the method of R. Brylinski and B. Kostant, we will
construct a representation
$\rho $ of
${\goth g}={\goth k}+{\goth  p}$ on the space of finite sums 
$${\cal O}(\Xi)_{\rm fin}=\sum _{m=0}^{\infty}{\cal O}_m(\Xi ),$$
such that, for all $X\in {\goth k}$, $\rho (X)=d\pi (X)$.
We define first a representation  $\rho$ of the subalgebra generated by $E,F, H$,
isomorphic to ${\goth sl}(2,{\bboard C})$. In particular
$$\rho (H)=d\pi (H)={d\over dt}\Big| _{t=0}\pi (\exp tH).$$
Hence, for $\phi\in {\cal O}_m(\Xi )$,
$\rho(H)\phi=({\cal E}-2m)\phi$,
where $\cal E$ is the Euler operator
$${\cal E}\phi (w,z)={d\over dt}\Big|_{t=0}\phi(w,e^tz).$$
One introduces two  operators $\cal M$ and $\cal D$. The operator $\cal M$ is a
multiplication operator:
$$({\cal M}\phi )(w,z)=w\phi (w,z),$$ 
which maps ${\cal O}_m(\Xi)$ into ${\cal O}_{m+1}(\Xi )$,
and $\cal D$ is  a differential operator: 
$$({\cal D}\phi)(w,z)={1\over w}\biggl(Q\Bigl({\partial\over \partial z}\Bigr)\phi\biggr)(w,z),$$
which maps 
${\cal O}_m(\Xi)$ into ${\cal O}_{m-1}(\Xi)$.
(Recall that ${\cal O}_{-1}(\Xi )=\{0\}$.) We denote by ${\cal M}^{\sigma}$ 
and ${\cal D}^{\sigma}$ the conjugate operators: 
$$
{\cal M}^{\sigma}=\pi(\sigma){\cal M}\pi(\sigma)^{-1}, \quad
{\cal D}^{\sigma}=\pi(\sigma){\cal D}\pi(\sigma)^{-1}.$$ 
Given a sequence $(\delta_m)_{m\in {\bboard N}}$ one defines  the diagonal operator 
$\delta$ on ${\cal O}(\Xi)_{fin}$   by  
$$\delta(\sum _m \phi_m)=\sum _m \delta_m\phi_m,$$
and  put 
$$\eqalign{
\rho(F)&={\cal M}-\delta\circ {\cal D}, \cr
\rho(E)&=\pi(\sigma)\rho_0(F)\pi(\sigma)^{-1}
={\cal M}^{\sigma}-\delta\circ {\cal D}^{\sigma}.\cr}$$
(Observe that, since ${\rm deg}\, Q=4$, then $Q$ is even, and $\sigma=\sigma^{-1}$.)

\tha Lemma 5.1|
$$\eqalign{
[\rho(H),\rho(E)]&=2\rho(E),\cr
[\rho(H),\rho(F)]&=-2\rho(F).\cr}$$
\fintha

\proof
Since
$$\eqalign{ 
\rho(H) {\cal M} & : \psi(z)w^m\mapsto ({\cal E}-2(m+1))\psi(z)w^{m+1},\cr
{\cal M} \rho(H) & : \psi(z)w^m \mapsto ({\cal E}-2m)\psi(z)w^{m+1},\cr}$$
one obtains $[\rho(H),{\cal M}]=-2{\cal M}$.
Since 
$$\eqalign{
\rho(H)\delta {\cal D} & : \psi(z)w^m\mapsto 
\delta_{m-1}({\cal E}-2(m-1))Q\Bigl({\partial\over\partial z}\Bigr)\psi(z)w^{m-1},\cr
\delta {\cal D}\rho(H) & : \psi(z)w^m \mapsto
\delta_{m-1}Q\Bigl({\partial\over \partial z}\Bigr)({\cal E}-2m)\psi(z)w^{m-1},\cr}$$
and, by using the identity 
$$[Q\Bigl({\partial\over \partial z}\Bigr),{\cal E}]=4Q\Bigl({\partial\over \partial z}\Bigr),$$
one gets 
$$[\rho (H),\delta {\cal D}]: \psi(z)w^m\mapsto
2\delta_{m-1}Q\Bigl({\partial \over \partial z}\Bigr)\psi(z)w^{m-1}.$$
Finally $[\rho(H),\rho(F)]=-2\rho(F)$.
Since  the operator $\delta$ commutes with $\pi(\sigma)$, and 
$\pi(\sigma)\rho(H)\pi(\sigma)^{-1}=-\rho(H)$,
we get also 
$[\rho(H),\rho(E)]=2\rho(E)$.
\hfill \qed

\bigskip

Let ${\bboard D}(V)^L$ denote the algebra of $L$-invariant differential operators on $V$. This algebra is commutative. If $V$ is simple and $Q=\Delta $, the determinant polynomial, then ${\bboard D}(V)^L$
is isomorphic to the algebra ${\cal P}({\bboard C}^r)^{{\goth S} _r}$ of symmetric polynomials in $r$ variables. The map
$$D\mapsto \gamma (D),\quad {\bboard D}(V)\to {\cal P}({\bboard C}^r)^{{\goth S}_r},$$
is the Harish-Chandra isomorphism (see Theorem XIV.1.7 in [Faraut-Kor\'anyi,1994]).
In general $V$ decomposes into simple ideals,
$$V=\bigoplus _{i=1}^s V_i,$$
and ${\bboard D}(V)^L$ is isomorphic to the algebra
$$\prod _{i=1}^s {\cal P}({\bboard C}^{r_i})^{{\goth S}_{r_i}}.$$
The isomorphism is given by
$$D\mapsto \gamma (D)=\bigl(\gamma _1(D),\ldots ,\gamma _s(D)\bigr),$$
where $\gamma _i$ is the isomorphism relative to the algebra $V_i$.
For $D\in {\bboard D}(V)^L$, we define the adjoint $D^*$ by
$D^*=J\circ D\circ J$, where $Jf(z)=f\circ j (z)=f(-z^{-1})$. Then
$\gamma (D^*)(\lambda )=\gamma (D)(-\lambda )$.
(See Proposition XIV.1.8 in [Faraut-Kor\'anyi,1994].)

In our setting we define the Maass operator ${\bf D}_{\alpha }$ as
$$D_{\alpha }=Q(z)^{1+\alpha }Q\Bigl({\partial \over \partial z}\Bigr)Q(z)^{-\alpha }.$$
It is $L$-invariant. We write
$$\gamma _{\alpha } (\lambda )=\gamma (D_{\alpha })(\lambda ).$$
If $V$ is simple and $Q=\Delta $, then
$$\gamma _{\alpha } (\lambda )=\prod _{i=1}^r \Bigl(\lambda _j-\alpha +\demi ({n\over r} -1)\Bigr),$$
([Faraut-Kor\'anyi,1994], p.296). If $V$ is simple and $Q=\Delta ^k$, then
$$\eqalign{
{\bf D}_{\alpha } &=\Delta ^{k+k\alpha }\Delta \Bigl({\partial \over \partial z}\Bigr)^k
\Delta (z)^{-k\alpha } \cr
&=\prod _{j=1}^k \Delta ^{k\alpha +k-j+1}\Delta \Bigl({\partial \over \partial z}\Bigr)
\Delta ^{-(k\alpha +k-j)},\cr}$$
and
$$\gamma _{\alpha } (\lambda )=\prod _{j=1}^r\big[\lambda _j-k\alpha +\demi ({n\over r}-1)\big]_k.$$
(We have used the Pochhammer symbol $[a]_k=a(a-1)\ldots (a-k+1)$.)

\th Proposition 5.2| 
In general
$$\gamma _{\alpha } (\lambda )=
\prod _{i=1}^s \prod _{j=1}^{r_i} \big[\lambda _j^{(i)}-k_i\alpha +\demi ({n_i\over r_i}-1)\big]_{k_i},$$
for $\lambda =(\lambda ^{(1)},\ldots ,\lambda ^{(s)})$, $\lambda ^{(i)}\in {\bboard C}^{r_i}$.
\finth

\vfill \eject

We say that the pair $(V,Q)$ has property (T) if there is a constant $\eta $ such that, for $X\in {\goth l}=Lie(L)$,
$$\Tr (X)=\eta  \tau (X).$$
In such a case, for $g\in L$,
$$\Det (g)=\gamma (g)^{\eta },$$
and, for $x\in V$,
$$\Det \bigl(P(x)\bigr)=Q(x)^{2\eta }.$$
Furthermore $Q(x)^{-\eta }m(dx)$ is an $L$-invariant measure on the symmetric cone
$\Omega \subset V_{\bboard R}$, and $H_0(z)=H(z)^{-2\eta }$.

Let $V=\oplus _{i=1}^s V_i$ be the decomposition of $V$ into simple ideals. Property (T)
is equivalent to the following: there is a constant $\eta $ such that
$${n_i\over r_i}=\eta  k_i\quad (i=1,\ldots ,s).$$
In fact, for $x\in V$, 
$$\Tr (T_x)=\sum _{i=1}^s {n_i\over r_i}\tr _i (x_i), \quad
\tau (T_x)=\sum _{i=1}^s k_i\tr _i(x_i),$$
with $x=(x_1,\ldots ,x_s)$, $x_i\in V_i$.

Property (T) is satisfied either if $V$ is simple, or if
$V={\bboard C}^p \oplus {\bboard C}^p$, and
$$Q(z)=(z_1^2+\cdots +z_p^2)(z_{p+1}^2+\cdots +z_{2p}^2).$$
Hence we get the following cases with property (T):

(1) $V={\bboard C}^n$, $Q(z)=(z_1^2+\cdots +z_n^2)^2$, and then
$${\goth g}={\goth sl}(n+2,{\bboard C}),\quad {\goth k}={\goth so}(n+2,{\bboard C}).$$

(2) $V={\bboard C}^p\oplus {\bboard C}^p$, and then
$${\goth g}={\goth so}(2p+4,{\bboard C}),\quad 
{\goth k}={\goth so}(p+2,{\bboard C})\oplus {\goth so}(p+2,{\bboard C}).$$

(3) $V$ is simple of rank 4, and $Q=\Delta $, the determinant polynomial. Then
$$({\goth g},{\goth k})=\bigl({\goth e}_6,{\goth sp}(8,{\bboard C})\bigr),\quad
\bigl({\goth e}_7,{\goth sl}(8,{\bboard C})\bigr),\quad
\bigl({\goth e}_8,{\goth so}(16,{\bboard C})\bigr).$$

\medskip

Observe that the case $V={\bboard C}^2$, $Q(z_1,z_2)=(z_1z_2)^2=z_1^2z_2^2$
belongs both to (1) and (2). This corresponds to the isomorphisms:
$${\goth sl}(4,{\bboard C})\simeq {\goth so}(6,{\bboard C}),\
{\goth so}(4,{\bboard C})\simeq {\goth so}(3,{\bboard C})\oplus {\goth so}(3,{\bboard C}).$$

\th Proposition 5.3|
The subspaces ${\cal O}_m(\Xi )$ are invariant under $[\rho (E),\rho (F)]$,
and the restriction of  $[\rho(E),\rho(F)]$ to ${\cal O}_m(\Xi )$ 
commutes with the $L$-action:
$$[\rho(E),\rho(F)]: {\cal O}_m(\Xi) \rightarrow  {\cal O}_m(\Xi), \quad
\psi(z)w^m\mapsto (P_m\psi)(z)w^m,$$
where $P_m$ is an $L$-invariant differential operator
on $V$ of degree $\leq 4$. It is given by 
$$P_m=\delta_m({\bf D}_{-1}-{\bf D}_{-m-1}^*)
+\delta_{m-1}({\bf D}_{-m}^*-{\bf D}_0).$$
\finth

\proof
Restricted to ${\cal O}_m(\Xi)$, 
$${\cal M}^{\sigma}{\cal D}={\bf D}_0,\quad
{\cal D}{\cal M}^{\sigma}={\bf  D}_{-1}, \quad
{\cal M}{\cal D}^{\sigma}={\bf D}_{-m}^*, \quad
{\cal D}^{\sigma}{\cal M}={\bf D}_{-m-1}^*.$$
It follows that the restriction of the operator $[\rho(E),\rho(F)]$ 
to ${\cal O}_m(\Xi)$ is given by 
$$\eqalignno{
[\rho(E),\rho(F)]&=[{\cal M}^{\sigma }
-\delta\circ{\cal D}^{\sigma},{\cal M}-\delta\circ {\cal D} ] \cr
&=[{\cal M},\delta \circ {\cal D}^{\sigma}]+[\delta \circ {\cal D},{\cal M}^{\sigma}] \cr
&={\cal M}\delta {\cal D}^{\sigma}-\delta {\cal D}^{\sigma}{\cal M}
+\delta {\cal D} {\cal M}^{\sigma}-{\cal M}^{\sigma} \delta \circ {\cal D} \cr
&=\delta_m({\cal D}{\cal M}^{\sigma}-{\cal D}^{\sigma}{\cal M})
+\delta_{m-1}({\cal M} {\cal D}^{\sigma}-{\cal M}^{\sigma}{\cal D}) \cr
&=\delta_m({\bf D}_{-1}-{\bf D}_{-m-1}^*)
+\delta_{m-1}({\bf D}_{-m}^*-{\bf D}_0).& \qed \cr}$$

\bigskip

By the Harish-Chandra isomorphism the operator $P_m$ corresponds to the polynomial
$p_m=\gamma (P_m)$,
$$p_m(\lambda )=\delta_m\bigl(\gamma_{-1}(\lambda)-\gamma_{-m-1}(-\lambda)\bigr)
+\delta_{m-1}\bigl(\gamma_{-m}(-\lambda)-\gamma_0(\lambda)\bigr).$$
The question is now whether it is possible to choose the sequence $(\delta_m)$
in such a way that $[\rho(E),\rho(F)]=\rho(H)$.  Recall that restricted
to ${\cal O}_m(\Xi)$, 
$$\rho(H)={\cal E}-2m,$$
where $\cal E$ is the Euler operator
$${\cal E}\phi(w,z)={d\over dt} \big|_{t=0}\phi(w,e^tz).$$
Then, by Proposition 5.3, it amounts to checking that, for every $m$,
$$p_m(\lambda)=\gamma ({\cal E})(\lambda )-2m.$$

\th Theorem 5.4|
It is possible to choose the sequence 
$(\delta_m)$ such that 
$$[\rho(H),\rho(E)]=2\rho(E), \quad [\rho(H),\rho(F)]=-2\rho(F),
\quad [\rho(E),\rho(F)]=\rho(H),$$
if and only if $(V,Q)$ has property (T), and then
$$\delta_m={A\over (m+\eta )(m+\eta +1)},$$  
where $A$ is a constant depending on  $(V,Q)$.
\finth

\proof 
a) Let us assume first that the Jordan algebra $V$ is simple of rank 4.
In such a case
$$\gamma _{\alpha } (\lambda )=\prod _{j=1}^4 \Bigl(\lambda _j-\alpha 
+\demi (\eta  -1)\Bigr)\quad (\eta ={n\over r})$$
(Proposition 5.2) . With $X_j=\lambda _j+\demi (\eta -1)$, the polynomial $p_m$ can be written
$$\eqalign{
p_m(\lambda )
&=\delta _m\Bigl(\prod _{j=1}^4 (X_j+1)-\prod _{j=1}^4 (X_j-m-\eta )\Bigr) \cr
&+\delta _{m-1}\Bigl(\prod _{j=1}^4 (X_j-m+1-\eta )-\prod _{j=1}^4X_j\Bigr).\cr}$$ 
Furthermore
$$\gamma ({\cal E})(\lambda )-2m=\sum _{j=1}^4 \lambda _j -2m 
=\sum _{j=1}^4X_j-2(m+\eta -1).$$

\th Lemma 5.5|
The identity in the four variables $X_j$
$$\eqalign{
&\alpha \Bigl(\prod _{j=1}^4 (X_j+1)-\prod _{j=1}^4 (X_j-b_j-1)\Bigr)
+\beta \Bigl(\prod _{j=1} ^4 (X_j-b_j)-\prod _{j=1}^4 X_j \Bigr)\cr
&=\sum _{j=1}^4 X_j+c \cr}$$
holds if and only if there is a constant $b$ such that
$$\eqalign{
&b_1=b_2=b_3=b_4=b,\ c=-2b, \cr
&\alpha ={1\over (b+1)(b+2)},\ \beta ={1\over b(b+1)}.\cr}$$
\finth

Hence we apply the lemma, and get $b=m+\eta -1$.
\hfill \qed

\bigskip

b) In the general case
$$\eqalign{
\gamma _{\alpha } (\lambda )
&=\prod _{i=1}^s \prod _{j=1}^{r_i} \big[\lambda _j^{(i)}-k_i\alpha +\demi ({n_i\over r_i}-1)\big]_{k_i}\cr
&=\prod _{i=1}^s\prod _{j=1}^{r_i}\prod _{k=1}^{k_i}
\Bigl(\lambda _j^{(i)}-k_i\alpha +\demi \bigl({n_i\over r_i}-1\Bigr)-(k-1)\Bigr) \cr
&=A\prod _{i=1}^s\prod _{j=1}^{r_i}\prod _{k=1}^{k_i}
\Bigl({\lambda _j^{(i)}\over k_i}-\alpha +{1\over 2k_i}\bigl({n_i\over r_i}-1\bigr)-{k-1\over k_i}\Bigr),\cr}$$
with $A=\prod _{i=1}^sk_i^{k_ir_i}$.
We introduce the notation
$$\eqalign{
X_{jk}^{(i)}&={\lambda _j^{(i)}\over k_i}+{1\over 2k_i}\bigl({n_i\over r_i}-1\bigr)-{k-1\over k_i},\cr
b_m^{(i)}&=m+{n_i\over k_ir_i}-1.\cr}$$
Then we obtain
$$\eqalign{
p_m(\lambda )
=&A \delta _m\Bigl(\prod _{i=1}^s\prod _{j=1}^{r_i}\prod _{k=1}^{k_i} (X_{jk}^{(i)}+1)
-\prod _{i=1}^s\prod _{j=1}^{r_i}\prod _{k=1}^{k_i} (X_{jk}^{(i)}-b_m^{(i)}-1) \Bigr)\cr
&+A\delta _{m-1}\Bigl(\prod _{i=1}^s\prod _{j=1}^{r_i}\prod _{k=1}^{k_i}(X_{jk}^{(i)}-b_m^{(i)})
-\prod _{i=1}^s\prod _{j=1}^{r_i}\prod _{k=1}^{k_i}(X_{jk}^{(i)})\Bigr), \cr}$$
and
$$\gamma ({\cal E})(\lambda )=\sum _{i=1}^s\sum _{j=1}^{r_i}\sum _{k=1}^{k_i}X_{jk}^{(i)}
-\demi \sum _{i=1}^s\sum _{j=1}^{r_i}\sum _{k=1}^{k_i}b_m^{(i)}.$$

If the rank of $V$ is equal to 4, then the $k_i$ are equal to 1, and the four variables
$X_{j1}^{(i)}$ are independant. 
By Lemma 5.5, Theorem 5.4 is proven in that case.

If the rank $r$ of $V$ is $<4$, then
$$X_{jk}^{(i)}=X_{j1}^{(i)}-{k-1\over k_i},$$
and there are only $r$ independant variables: $X_{j1}^{(i)}$.
In that case Theorem 5.4 is proven by using an alternative form of Lemma 5.5:
\hfill \qed

\th Lemma 5.6|
To a partition $k=(k_1,\ldots ,k_{\ell })$ of 4 and length $\ell $:
$$k_1+\dots +k_{\ell }=4,$$
and the numbers $\gamma _{ij}$ ($1\leq i\leq \ell $, $1\leq j\leq k_i-1$), 
one associates
the polynomial $F$ in the $\ell $ variables $T_1,\ldots ,T_{\ell }$:
$$F(T_1,\ldots ,T_{\ell })=\prod _{i=1}^{\ell } T_i\prod _{j=1}^{k_i-1}(T_i+\gamma _{ij}).$$
Given $\alpha ,\beta,c\in {\bboard R}$, and $b_1,\ldots b_{\ell }\in {\bboard R}$, then
$$\eqalign{
&\alpha \bigl(F(T_1+1,\ldots ,T_{\ell }+1) 
-F(T_1-b_1-1,\ldots ,T_{\ell }-b_{\ell }-1)\bigr)\cr 
&+\beta \bigl(F(T_1-b_1,\ldots ,T_{\ell }-b_{\ell })-F(T_1,\ldots ,T_{\ell }\bigr) 
=\sum _{i=1}^{\ell }T_i+c \cr}$$
is an identity in the variables $T_1,\ldots ,T_{\ell }$ if and only if
there exists $b$ such that 
$$b_1=\ldots =b_{\ell }=b,\ \alpha ={1\over (b+1)(b+2)},\ \beta ={1\over b(b+1)},$$
and
$$c=\sum _{i=1}^{\ell }\sum _{j=1}^{k_i-1}\gamma _{ij}-2b.$$
\finth

\bigskip
 
For $p\in {\goth p}$, define the multiplication operator ${\cal M}(p)$ given by
$$\bigl({\cal M}(p)\phi \bigr)(w,z)=wp(z)\phi (w,z).$$
Observe that ${\cal M}(1)={\cal M}$. Then, for $g\in K$,
$${\cal M}\bigl(\kappa (g)p\bigr)=\pi (g){\cal M}(p)\pi (g^{-1}).$$
In fact
$$\bigl({\cal M}(p)\pi (g^{-1})\phi \bigr) (w,z)
=wp(z)\phi \bigl(\mu (g,z)w,g\cdot z),$$
and
$$\eqalign{
&\bigl( \pi (g){\cal M}(p)\pi (g^{-1})\phi \bigr)(w,z)\cr
&=\mu (g^{-1},z)wp(g^{-1}\cdot z)\phi \bigl( \mu (g^{-1},z)\mu (g,g^{-1}\cdot z)w,g^{-1}g\cdot z\bigr) \cr
&=w\bigl(\kappa (z)p\bigr)(z)\phi (w,z)={\cal M}\bigl(\kappa (g)p\bigr)\phi (w,z).\cr}$$

\th Proposition 5.7|
There is a unique map 
$${\goth p}\to {\rm End}\bigl({\cal O}_{\rm fin }(\Xi )\bigr),\quad p\mapsto {\cal D}(p),$$
such that ${\cal D}(1)={\cal D}$, and, for $g\in K$,
$${\cal D}\bigl(\kappa (g)p\bigr)=\pi (g){\cal D}(p)\pi (g^{-1}).$$
\finth 

\proof
Recall that, for $g\in P_{\rm max}$,
$$\bigl(\kappa (g)p\bigr)(z)=\chi (g)p(g^{-1}\cdot z),$$
and
$$\bigl(\pi (g)\phi \bigr)(w,z)=\phi \bigl(\chi (g)w,g^{-1}\cdot g\bigr).$$
Let us show that, for $g\in P_{\rm max}$,
$$\pi (g){\cal D}\pi (g^{-1})=\chi (g){\cal D}.$$
Observe first that, for $\ell \in L$ and a smooth function $\psi $ on $V$,
$$Q\Bigl({\partial \over \partial z}\Bigr)\bigl(\psi (\ell \cdot z)\bigr)=\gamma (\ell )
\Bigl(Q\bigl({\partial \over \partial z}\bigr)\psi \Bigr)(\ell \cdot z).$$
Therefore, for $g\in P_{\rm max}$,
$$\eqalign{
{\cal D}\pi (g^{-1})\phi (w,z)
&={1\over w}Q\Bigl({\partial \over \partial z}
\Bigl(\phi \bigl(\chi (g^{-1})w,g\cdot z\bigr)\Bigr) \cr
&={1\over w}\chi (g)^2\Bigl(Q\bigl({\partial \over \partial z}\bigr)\phi \Bigr)\bigl(\chi (g^{-1})w,g\cdot z\bigr),\cr}$$
and
$$\bigl(\pi (g){\cal D}\pi (g^{-1})\phi \bigr)(w,z)
={1\over \chi (g)w}\chi (g)^2\Bigl(Q\bigl({\partial \over \partial z}\bigr)\phi \Bigr)(w,z)
=\chi (g){\cal D}\phi (w,z).$$
It follows that the vector subspace in ${\rm End}\bigl({\cal O}_{\rm fin}(\Xi )\bigr)$
generated by the endomorphisms $\pi (g){\cal D}\pi (g^{-1})$ ($g\in K$) is a representation space
for $K$ equivalent to $\goth p$. (See Theorem 3.10 in [Brylinski-Kostant,1994].)
Hence there exists a unique $K$-equivariant map $p\mapsto {\cal D}(p)$ such that
${\cal D}(1)={\cal D}$.

\hfill \qed

\bigskip

For $p\in {\goth p}$, define 
$$\rho (p)={\cal M}(p)-\delta {\cal D}(p).$$
Observe that this definition is consistent with the definition of $\rho (E)$ and $\rho (F)$.
Recall that, for $X\in {\goth k}$, $\rho (X)=d\pi (X)$. Hence we get a map
$$\rho :{\goth g}={\goth k}\oplus {\goth p}\to {\rm End}\bigl({\cal O}(\Xi )_{\rm fin}\bigr).$$

\th Theorem 5.8|
Assume that Property $(T)$ holds. Fix $(\delta _m)$ as in Theorem 5.4.

{\rm (i)} $\rho $ is a representation of the Lie algebra $\goth g$ on ${\cal O}(\Xi )_{\rm fin}$.

{\rm (ii)} The representation $\rho $ is irreducible.
\finth

\proof
(i) Since $\pi $ is a representation of $K$, for $X,X'\in {\goth k}$,
$$[\rho (X),\rho (X')]=\rho ([X,X']).$$
It follows from Proposition 5.7 that, for $X\in {\goth k},p\in {\goth p}$,
$$[\rho (X),\rho (p)]=\rho ([X,p]).$$
It remains to show that, for $p,p'\in {\goth p}$,
$$[\rho (p),\rho (p']]=\rho ([p,p'].$$
By Theorem 5.4, $[\rho (E),\rho (F)]=\rho (H)$. Then this follows from Lemma 3.6 in [Brylinski-Kostant,1995]: consider the map
$$\tau :{\bigwedge }^2{\goth p}\to {\rm End}\bigl({\cal O}(\Xi )_{\rm fin},$$
defined by
$$\tau (p\wedge p')=[\rho (p),\rho (p')]-\rho ([p,p']).$$
We know that $\tau (E\wedge F)=0$. It follows that, for $g\in K$,
$$\tau \bigl(\kappa (g)E\wedge \kappa (g)F)=0.$$
Since the representation $\kappa $ is irreducible, and $E$ and $F$ are highest and lowest vectors with respect to $P$, the vector $E\wedge F$ is cyclic in ${\bigwedge }^2{\goth p}$ for the action of $K$. Therefore $\tau \equiv 0$.

(ii) Let ${\cal V}\ne\{0\}$ be a  $\rho({\goth g})$-invariant subspace of 
${\cal O}(\Xi)_{\rm fin}$. Then  $\cal V$ is $\rho({\goth k})$-invariant. 
As ${\cal O}(\Xi)_{\rm fin}=\sum _{m=0}^{\infty }{\cal O}_m(\Xi)$ and as the 
subspaces ${\cal O}_m(\Xi)$ are  $\rho({\goth k})$-irreducible, then there exists 
${\cal I}\subset {\bboard N}$ (${\cal I}\not=\emptyset $) such that  
${\cal V}=\sum _{m\in {\cal I}}{\cal O}_m(\Xi)$.
Observe that if ${\cal V}$ contains ${\cal O}_m(\Xi)$, 
then it contains ${\cal O}_{m+1}(\Xi)$ too. 
In fact denote by $\phi _m$ the function in ${\cal O}_m(\Xi )$ defined by
$\phi _m(w,z)=w^m$. As ${\cal D} \phi_{m}=0$, it follows that  
$$\rho(F)\phi_{m}={\cal M}\phi_{m}=\phi_{m+1},$$
and $\rho(F)\phi_{m}$ belongs to  ${\cal O}_{m+1}(\Xi)$, therefore
${\cal O}_{m+1}(\Xi) \subset {\cal V}$. 
Denote by $m_0$ the minimum of the $m$ such that ${\cal O}_m(\Xi)\subset {\cal V}$, then
$${\cal V}=\bigoplus _{m=m_0}^{\infty}{\cal O}_m(\Xi).$$
The function $\phi(w,z)=Q(z)^mw^m$ belongs to ${\cal O}_m(\Xi)$, and
$$\rho(F)\phi(w,z)=Q(z)^mw^{m+1}
-\delta_{m-1}Q\Bigl({\partial \over \partial z}\Bigr)Q(z)^mw^{m-1}.$$
By the Bernstein identity (Proposition 3.1)
$$Q\Bigl({\partial \over \partial z}\Bigr)Q(z)^m=B(m)Q(z)^{m-1},$$
and since $B(m)>0$ for $m>0$, it follows that, if ${\cal O}_m(\Xi) \subset {\cal V}$
with $m>0$, then ${\cal O}_{m-1}(\Xi) \subset {\cal V}$. 
Therefore $m_0=0$ and ${\cal V}={\cal O}(\Xi)_{\rm fin}$.
\hfill \qed

\section 6. The unitary representation of the Kantor-Koecher-Tits \break group|
We consider, for a sequence $(c_m)$ of positive numbers, an inner product on 
${\cal O}(\Xi)_{\rm fin}$  such that 
$$\Vert\phi\Vert^2=\sum\limits_{m=0}^{\infty}{1\over c_m}\Vert
\psi_m\Vert_m^2,$$
for
$$\phi(w,z)=\sum_{m=0}^{\infty}\psi_m(z)w^m.$$
This inner product is invariant under $K_{\bboard R}$. 
We assume that Property (T) holds,
and we will determine the sequence
$(c_m)$ such that this inner product  is invariant under  the representation $\rho $ restricted to ${\goth g}_{\bboard R}$. We denote
by $\cal H$ the Hilbert space completion of ${\cal O}(\Xi)_{\rm fin}$ with
respect to this inner product. We will assume $c_0=1$.

The Bernstein polynomial $B$ is of degree 4, and vanishes at 0 and $\alpha _1=1-\eta $. 
Let $\alpha _2$ and $\alpha _3$ be the two remaining roots:
$$B(\alpha )=A\alpha (\alpha -\alpha _1)(\alpha -\alpha _2)(\alpha -\alpha _3).$$

\bigskip

(1) $V={\bboard C}^n$, $Q(z)=(z_1^2+\cdots +z_n^2)^2$. Then
$$B(\alpha )=A\alpha \Bigl(\alpha -{1\over 2}\Bigr)\Bigl(\alpha +{n-4\over 4}\Bigr)
\Bigl(\alpha + {n-2\over 4}\Bigr).$$
$A=2^4$ if $n\geq 2$, $A=4^4$ if $n=1$.

\medskip

(2) $V=(z_1^2+\cdots +z_p^2)(z_{p+1}^2+\cdots +z_{2p}^2)$. Then
$$B(\alpha )=\alpha ^2\Bigl(\alpha +{p-2\over 2}\Bigr)^2.$$

\medskip

(3) $V$ is simple of rank 4, complexification
of $V_{\bboard R}=Herm(4,{\bboard F})$, $Q(z)=\Delta (z)$, the determinant polynomial. Then
$$B(\alpha )=\alpha \Bigl(\alpha +{d\over 2}\Bigr)\Bigl(\alpha +2{d\over 2}\Bigr)
\Bigl(\alpha +3{d\over 2}\Bigr),$$
where $d={\rm dim }_{\bboard R}{\bboard F}$.

\bigskip

Here are the non zero roots of the Bernstein polynomial: 

\def\tvi{\vrule height 12pt depth 1pt width 0pt}

\def\traithorizontal{\noalign{\hrule}}

$$\vbox{\offinterlineskip \halign{
\tvi # & \quad # & \quad # & \quad #  & \quad # \cr
& $\eta $ & $\alpha _1 $ & $\alpha _2$ & $\alpha _3$ \cr
\traithorizontal 
(1) & ${n\over 4}$ & $-{n-4\over 4}$ & ${1\over 2}$ & $-{n-2\over 4}$ \cr
(2) & ${p\over 2}$ & $-{p-2\over 2}$ & $0$ & $-{p-2\over 2}$ \cr
(3) & $1+3{d\over 2}$ & $-3{d\over 2}$ & $-{d\over 2}$ & $-2{d\over 2}$ \cr
}}$$

\th Theorem 6.1|
{\rm (i)} The inner product of $\cal H$ is ${\goth g}_{\bboard R}$-invariant 
if 
$$c_m={(\eta +1)_m\over (\eta +\alpha _2 )_m(\eta +\alpha _3)_m}{1\over m!}.$$

{\rm (ii)} The  reproducing kernel of $\cal H$ is given by 
$${\cal K}(\xi,\xi')={}_1F_2(\eta +1;\eta +\alpha _2,\eta +\alpha _3;H(z,z')w\overline{w'}),$$
for $\xi=(w,z)$, $\xi'=(w',z')$.
 \finth

\proof
(i) Recall that
$${\goth p}_{\bboard R}=\{p\in {\goth p}\mid \beta (p)=p\},$$
where $\beta $ is the conjugation of $\goth p$, we introduced at the end of Section 4.
Recall also that 
$$\beta (\kappa (g)p)=\kappa \bigl(\alpha (g)\bigr)\beta (p).$$
The inner product of $\cal H$ is ${\goth g}_{\bboard R}$-invariant if and only if, for every $p\in {\goth p}$,
$$\rho (p)^*=-\rho \bigl(\beta (p)\bigr).$$
But this is equivalent to the single condition
$$\rho (E)^*=-\rho (F).$$
In fact, assume that this condition is satisfied. Then, for $p=\kappa (g)E$, ($g\in K$),
$$\rho (p)=\pi (g)\rho (E)\pi (g^{-1}),\quad \rho (p)^*=-\pi (g^{-1})^*\rho (F)\pi (g)^*.$$
Since $\pi (g)^*=\pi \bigl(\alpha (g)\bigr)^{-1}$, we get
$$\eqalign{
\rho (p)^*
&=-\pi \bigl(\alpha (g)\bigr)\rho (F)\pi \bigl(\alpha (g^{-1})\bigr)
=-\rho \bigl(\kappa (\alpha (g))F\bigr) \cr
&=-\rho \bigl(\kappa (\alpha (g))\beta (E)\bigr) 
=-\rho \bigl(\beta (\kappa (g)E)\bigr)=-\rho \bigl(\beta (p)\bigr).\cr}$$
Finally observe that the vector $E$ is cyclic in $\goth p$ for the $K$-action.

The condition $\rho (E)^*=-\rho (F)$ is equivalent to: for $m\geq 0$,
$\phi\in {\cal O}_{m+1}(\Xi), \phi'\in {\cal O}_m(\Xi)$,
$${1\over c_{m+1}}(\phi\mid{\cal M}^{\sigma}\phi')_{m+1}
={1\over c_m}\delta_m({\cal D}\phi\mid\phi')_m.$$
Recall that $m_0(dz)=H_0(z)m(dz)$ with
$$H_0(z)=H(z)^{-2\eta },$$
and the norm of $\tilde{\cal O}_m(V)$ can be written
$$\Vert \psi\Vert_m^2={1\over a_m}\int_V\vert\psi(z)\vert^2H(z)^{-m-2\eta }m(dz).$$
Then, the  required condition of invariance becomes 
$$\eqalign{
&{1\over c_{m+1}a_{m+1}}\int_V\psi(z)\overline{Q(z)}
\overline{\psi'(z)}H(z)^{-(m+1)-2\eta }m(dz) \cr
&={\delta_m\over c_ma_m}
\int_V (Q\Bigl({\partial \over \partial z}\Bigr)\psi)(z)\overline{\psi'(z)}H(z)^{-m-2\eta }m(dz).\cr}$$
By  integrating by parts:
$$\eqalign{
&\int_V(Q\Bigl({\partial\over \partial z}\Bigr)\psi)(z)\overline
{\psi'(z)}H(z)^{-m-2\eta }m(dz) \cr
&=\int_V\psi(z)\overline{\psi'(z)}
\biggl(Q\Bigl({\partial\over \partial z}\Bigr)H(z)^{-m-2\eta }\biggr)m(dz),\cr}$$
and, by the relation
$$Q\Bigl({\partial \over \partial z}\Bigr) H(z)^{-m-2\eta } 
=B(-m-2\eta )\overline{Q(z)}H(z)^{-(m+1)-2\eta },$$
the condition can be written
$${1\over c_{m+1}}={a_{m+1}\over a_m}\delta _mB(-m-2\eta ){1\over c_m}.$$
From Proposition 3.2 it follows that
$${a_{m+1}\over a_m}={B(-m-\eta )\over B(-m-2\eta )}.$$
We obtain finally
$${c_{m+1}\over c_m}={m+\eta +1\over (m+\eta +\alpha _2)(m+\eta +\alpha _3)(m+1)},$$
and, since $c_0=1$,
$$c_m={(\eta +1)_m\over (\eta +\alpha _2 )_m(\eta +\alpha _3)_m}{1\over m!}.$$

(ii) By Theorem 2.5 the reproducing kernel of $\cal H$ is given by
$$\eqalign{
{\cal K}(\xi,\xi')
&=\sum _{m=0}^{\infty }
c_m H(z,z')^{m}w^m\overline{w'}^m \cr
&={}_1F_2\bigl(\eta +\alpha _2,\eta +\alpha _3;\eta +1;H(z,z')w\overline{w'}\bigr),\cr}$$
with $\xi=(w,z)$, $\xi'=(w',z')$.
\hfill \qed 

\bigskip

We will see that the Hilbert space $\cal H$ is a pseudo-weighted Bergman space. By this we mean 
that the norm is given by an integral of $|\phi |^2$ with respect to a weight taking both positive and negative values.
The weight involves a Meijer $G$-function:
$$G(u)={1\over 2i\pi }\int _{c-i\infty }^{c+i\infty }
{\Gamma (\beta _1+s)\Gamma (\beta _2+s)\Gamma (\beta _3+s)\over \Gamma (\alpha +s)}u^{-s}ds,$$
where $\alpha ,\beta _1,\beta _2,\beta _3$ are real numbers,
and $c>\sigma =-\inf \{\beta _1,\beta _2,\beta _3\}$.
This function is denoted by
$$G(u)=G_{1,3}^{3,0}\Bigl(x\big|\matrix{\alpha & & \cr \beta _1 & \beta _2 & \beta _3 \cr}\Bigr)$$
(see for instance [Mathai,1993]). By the inversion formula for the Mellin transform
$$\int _0^{\infty } G(u)u^{s-1}du=
{\Gamma (\beta _1+s)\Gamma (\beta _2+s)\Gamma (\beta _3+s)
\over \Gamma (\alpha +s)},$$
for ${\rm Re}\, s>\sigma$, and the integral is absolutely convergent.
If the numbers $\beta _1,\beta _2,\beta _3$ are distinct, then
$$G(u)=\varphi _1(u)u^{\beta _1}+\varphi _2(u)u^{\beta _2}+\varphi _3(u)u^{\beta _3},$$
where $\varphi _1,\varphi _2,\varphi _3$ are holomorphic near 0.
($\varphi _1,\varphi _2,\varphi _3$ are $_1F_2$ hypergeometric functions.)

\medskip

The function $G$ may be not positive on $]0,\infty [$, but is positive for $u$ large enough.
In fact
$$G(u)\sim \sqrt{\pi }u^{\theta }e^{-2\sqrt{u}}\quad (u\to \infty ),$$
where
$$\theta =\beta _1+\beta _2+\beta _3-\alpha -\demi .$$
([Paris-Wood,1986], Theorem 3, p.32.)

Now take 
$$\alpha =\eta -1,\ \beta _1=2\eta -1,\ \beta _2 =2\eta +a-1,\ \beta _3=2\eta +b-1.$$

\bigskip

$$\vbox{\offinterlineskip \halign{
\tvi # & \quad # & \quad # & \quad #  & \quad # \cr
& $\alpha $ & $\beta _1 $ & $\beta _2$ & $\beta _3$ \cr
\traithorizontal 
(1) & ${n\over 4}-1$ & ${n-2\over 2}$ & ${n-1\over 2}$ & ${n-2\over 4}$ \cr
(2) & ${p\over 2}-1$ & $p-1$ & $p-1$ & ${p\over 2}$ \cr
(3) & $3{d\over 2}$ & $3d+1$ & $5{d\over 2}+1$ & $2d+1$ \cr
}}$$

\bigskip

The Mellin transform of $G$ vanishes at $-\alpha $, with changing sign. 
One can check that
$-\alpha >\sigma $ in all cases.
Therefore there are real values $s>\sigma $ for which the integral
$$\int _0^{\infty } G(u)u^{s-1}du<0.$$
This implies that the function $G$ takes negative values on $]0,\infty [$.

\th Theorem 6.2|
For $\phi \in {\cal H}$,
$$\|\phi \|^2=\int _{{\bboard C}\times V}
|\phi (w,z)|^2p(z,w)m(dw)m_0(dz),$$
with
$$p(w,z)=CG\bigl(|w|^2H(z)\bigr) H(z).$$
The integral is absolutely convergent.
\finth

\proof
We will follow the proof of Theorem 5.7 in [Brylinski,1997].

a) From the proof of Theorem 6.1 it follows that
$$\eqalign{
{1\over a_mc_m}
&={(2\eta )_m(2\eta +\alpha _2)_m(2\eta +\alpha _3)_m\over (\eta )_m}\cr
&=C{\Gamma (2\eta +m)\Gamma (2\eta +\alpha _2+m)\Gamma (2\eta +\alpha _3+m)
\over \Gamma (\eta +m)} \cr
&=C\int _0^{\infty } G(u)u^mdu.\cr}$$
(One checks that $\sigma <1$, {\it i.e.} $G$ is integrable.) By the computation we did for the proof of Theorem 2.6, we obtain, 
for $\phi (w,z)=w^m\psi (z)\in {\cal O}_m$,
$$\int _{{\bboard C}\times V}|\phi (w,z)|^2p(z,w)m(dw)m_0(dz)=\|\phi \|^2.$$
Furthermore, if $\phi \in {\cal O}_m$, $\phi '\in {\cal O}_{m'}$, with $m\not=m'$,
$$\int _{{\bboard C}\times V}\phi (w,z)\overline{\phi '(w,z)}m(dw)m_0(dz)=0.$$
It follows that, for $\phi \in {\cal O}_{\rm fin}$,
$$\int _{{\bboard C}\times V}|\phi (w,z)|^2p(z,w)m(dw)m_0(dz)=\|\phi \|^2.$$
The computation is justified by the fact that, for $s>\sigma $,
$$\int _0^{\infty } |G(u)|u^{s-1}du<\infty .$$

b) Let us consider the weighted Bergman space ${\cal H}^1$ whose norm is given by
$$\|\phi \|_1^2=\int _{{\bboard C}\times V}|\phi (w,z)|^2|p(w,z)|m(dw)m_0(dz).$$
By Theorem 2.6,
$$\|\phi \|_1^2=\sum _{m=0}^{\infty } {1\over c_m^1}\|\psi _m\|_m^2,$$
with 
$${1\over a_mc_m^1}=C \int _0^{\infty }|G(u)|u^mdu.$$
Obviously $c_m^1\leq c_m$, therefore ${\cal H}^1\subset {\cal H}$.
We will show that ${\cal H}\subset {\cal H}^1$. For that we will prove that there is a constant $A$ such that
$$c_m\leq A c_m^1.$$
As observed above there is $u_0\geq 0$ such that $G(u)\geq 0$, for $u\geq u_0$, and then
$$\int _0^{\infty } |G(u)|u^m
\leq \int _0^{\infty } G(u)u^mdu+2\int _0^{u_0}|G(u)|u^mdu.$$
Hence
$${1\over c_m^1}\leq {1\over c_m}+2a_mu_0^m\int _0^{u_0}|G(u)|du.$$
By the formula we gave at the beginning of a), the sequence $a_mc_mu_0^m$ is bounded. Therefore there is a constant $A$ such that
$${1\over c_m^1}\leq A {1\over c_m},$$
and this implies that ${\cal H}\subset {\cal H}_1$.
\hfill \qed

\vfill \eject

Let $\tilde{G_{\bboard R}}$ be the connected and simply
connected Lie group with Lie algebra ${\goth g}_{\bboard R}$ and denote by $\tilde
K_{\bboard R}$ the subgroup of $\tilde G_{\bboard R}$ with Lie algebra 
${\goth k}_{\bboard R}$. 
It is a covering of $K_{\bboard R}$. We denote by $s: \tilde K_{\bboard R}
\rightarrow K_{\bboard R}, g\mapsto s(g)$ the canonical surjection. 

\th Theorem 6.3| 
{\rm (i)} There is a unique unitary irreducible
representation $\tilde\pi$  of $\tilde G_{\bboard R}$ on $\cal H$ such that
d$\tilde\pi=\rho $.
And, for all $k\in \tilde K_{\bboard R}$, $\tilde\pi(k)=\pi(s(k))$.

{\rm (ii)} The representation $\tilde\pi$ is spherical. 
\finth

\proof
(i) Notice that if the operators $\rho(E+F)$ and
$\rho(i(E-F))$ are skew-symmetric, then for each $p\in {\goth p}_{\bboard R}$, the
operator $\rho (p)$ is skew-symmetric. In fact, since the ${\goth sl}_2$-triple
$(E,F,H)$ is strictly normal (see [Sekiguchi,1987]), which means that $H\in i{\goth k}_{\bboard R}, E+F\in 
{\goth p}_{\bboard R}, i(E-F) \in {\goth p}_{\bboard R}$,   and  since 
${\goth  p}={\cal U}({\goth k})E$,  hence 
${\goth p}_{\bboard R}={\cal U}({\goth k}_{\bboard R})(E+F)+{\cal U}({\goth k}_{\bboard R})(i(E-F))$,  and the assertion follows.

Now,   by Nelson's criterion, it is enough to prove that the
operator $\rho({\cal L})$ is essentially self-adjoint where $\cal L$ is the
Laplacian of ${\goth g}_{\bboard R}$. Let's consider a basis
$\{X_1,\ldots,X_k\}$ of ${\goth k}_{\bboard R}$  and a basis
$\{p_1,\ldots,p_l\}$ of ${\goth p}_{\bboard R}$, orthogonal with respect to the Killing form. 
As ${\goth g}_{\bboard R}={\goth k}_{\bboard R}+{\goth p}_{\bboard R}$ 
is the Cartan decomposition of ${\goth g}_{\bboard R}$, then the
Laplacian and the Casimir operators of ${\goth g}_{\bboard R}$ are given by 
$$\eqalign{
{\cal L}&=X_1^2+\ldots+X_k^2+p_1^2+\ldots+p_l^2,\cr
{\cal C}&=X_1^2+\ldots+X_k^2-p_1^2-\ldots-p_l^2.\cr}$$
It follows that ${\cal L}=2(X_1^2+\ldots+X_k^2)-{\cal C}$ and
$\rho({\cal L})=2\rho(X_1^2+\ldots+X_k^2)-\rho({\cal C})$.  
Since $\rho(X_1^2+\ldots+X_k^2)=d\pi(X_1^2+\ldots+X_k^2)$ and as
$\pi$ is a  unitary  representation of $K_{\bboard R}$,  hence the image
$\rho(X_1^2+\ldots+X_k^2)$ of the Laplacian of ${\goth k}_{\bboard R}$ is essentially
self-adjoint. Moreover, since the dimension of  ${\cal O}(\Xi)_{\rm fin}$
 is countable, then the commutant of
$\rho$, which is a division algebra over $\bboard C$,  has a countable
dimension too, and  is equal to $\bboard C$ (see [Cartier,1979], p.118). It follows that
$\rho({\cal C})$ is scalar.  We deduce that   $\rho({\cal L})$ is essentially
self-adjoint and that the irreducible representation 
$\rho$ of ${\goth g}_{\bboard R}$ integrates to an irreducible unitary representation of
$\tilde G_{\bboard R}$, on the Hilbert space $\cal H$.

(ii) The space ${\cal O}_0(\Xi)$ reduces to the constant
functions which are the $K$-fixed vectors.
\hfill \qed

\bigskip

We don't know whether the representation $\tilde \pi $ goes down to a representation of a real Lie group $G_{\bboard R}$ with $K_{\bboard R}$ as a maximal compact subgroup.

\vfill \eject

\def\article#1|#2|#3|#4|#5|#6|
       {{\noindent \hangindent=5mm \hangafter=1
      {\petcap #1} (#2). {\rm #3}, {\sl #4}, {\bf #5}, #6.\par}\medskip}

\def\livre#1|#2|#3|#4|
     {{\noindent \hangindent=5mm\hangafter=1
     {\petcap #1} (#2). {\rm #3}. {\sl #4}.\par}\medskip}

\centerline{\bf References}
\vskip 6pt
\rm

\article
D. Achab|2000|Alg\`ebres de Jordan de rang 4 et repr\'esentations
minimales|Advances in Mathematics|153|155-183|

\article
D. Achab|2011|Construction process for simple Lie algebras|
Journal of Algebra|325|186-204| 

\article
B. Allison|1979|Models of isotropic simple Lie algebras|Comm. in Alg.|7|1835-1875|

\article
B. Allison|1990|Simple structurable algebras of skew dimension one|
Comm. in Alg.|18|1245-1279|

\article
B. Allison and J. Faulkner|1984|A Cayley-Dickson Process for a class
of  structurable algebras|Trans. Amer. Math. Soc.|283|185--210|

\article
R. Brylinski|1997|Quantization of the 
4-dimensional nilpotent orbit of $SL(3,{\bboard R})$|
Canad. J. Math.|49|916-943|

\article
R. Brylinski|1998|Geometric quantization of real minimal nilpotent orbits.Symplectic geometry| 
Differential Geom. Appl.|9|5-58|

\article
R. Brylinski and B. Kostant|1994|
Minimal representations, geometric quantization, and unitarity|
Proc. Nat. Acad. USA|91|6026-6029|

\livre
R. Brylinski and B. Kostant|1995|
Lagrangian models of minimal representations of $E_8$, $E_7$ and $E_8$ 
in Functional Analysis on the Eve of the 21st Century. In honor of I.M. Gelfand's 80th Birthday, 13-53,
Progress in Math.131|Birkh\"auser|

\livre
R. Brylinski and B. Kostant|1997|Geometric quantization and holomorphic half-form models of unitary minimal representations I, II |Preprint|

\article
P. Cartier|1979|Representations of $p$-adic groups in  {\rm Automorphic
forms, representations and $L$-functions}|Proc. Symposia in Pure Math.|31{\rm.1}|111-155|

\livre
J. Faraut and A. Kor\'anyi|1994|Analysis on symmetric cones| Oxford University Press|

\livre
J. Faraut and S. Gindikin|1996|Pseudo-Hermitian symmetric spaces of
tube type, in  
{\rm Topics in Geometry (S. Gindikin ed.)}. Progress in non linear differential
equations and their applications, {\bf 20}, 123-154|Birkh\"auser|

\article
R. Goodman|2008|Harmonic analysis on compact symmetric spaces : the
legacy of Elie Cartan and Hermann Weyl in Groups and analysis|
London Math. Soc. Lecture Note|354|1-23|

\article
T. Kobayashi and G. Mano|2007|Integral formula of the unitary inversion operator for the minimal
representation of $O(p,q)$|
Proc. Japan Acad. Ser. A Math. Sci.|83|27--31|

\livre
T. Kobayashi and G. Mano|2008|The Schr\"odinger model for the minimal
representation of the indefinite orthogonal group $O(p,q)$|University of Tokyo, Graduate
School of Mathematical Sciences. Preprint,
to appear in {\it Memoirs of Amer. Math. Soc.}|

\article
T. Kobayashi and B. \O rsted|2003|Analysis on the minimal representation of $O(p,q)$.
I. Realization via conformal geometry|Adv. Math.|180|486--512|

\livre
A.M. Mathai|1993|A Handbook of Generalized Special Functions for
Statistical and Physical Sciences|Oxford University Press|

\article
K. McCrimmon|1978|Jordan algebras and their applications|
Bull. A.M.S.|84|612-627|

\article
M. Pevzner|2002|Analyse conforme sur les alg\`ebres de Jordan|    
J. Austral. Math. Soc.|73|1-21| 

\livre
R.B. Paris and A.D. Wood|1986|Asymptotics of high order differential
equations|Pitman Research Notes in Math Series, vol.
\bf 129, \rm Longman Scientific and Technical|Harlow|

\article
J. Rawnsley and S. Sternberg|1982|On representations associated to the minimal nilpotent coadjoint orbit of $SL(3,{\bboard R})$|Amer. J. Math.|104|1153--1180|

\article
J Sekiguchi|1987|Remarks on nilpotent orbits of a symmetric pair|
Jour. Math. Soc. Japan|39|127--138|

\article
P. Torasso|1983|Quantification g\'eom\'etrique et representations de $SL_3({\bboard R})$|
Acta Mathematica|150|153--242|

\vskip 1cm

\noindent
Institut de Math\'ematiques de Jussieu

\noindent
Universit\'e Pierre et Marie Curie

\noindent
4 place Jussieu, case 247, 75252 Paris cedex 05

\noindent
{\tt achab@math.jussieu.fr,\quad faraut@math.jussieu.fr}

\end